\documentclass{article}
\usepackage[margin=1.3in, top=1.3in, bottom=1.3in]{geometry}
\usepackage{graphicx}
\usepackage{enumitem}
\usepackage{amssymb}
\usepackage{amsthm}
\usepackage{amsmath}
\usepackage{latexsym}
\usepackage{color}
\usepackage[hypertexnames=false]{hyperref}
\usepackage{appendix}

\usepackage[normalem]{ulem}

\newcommand{\st}[1]{\textcolor{red}{\ifmmode\text{\sout{\ensuremath{#1}}}\else\sout{#1}\fi}}

\newcommand*{\Nwarrow}{\rotatebox[origin=c]{-45}{\(\Longrightarrow\)}}
\newcommand*{\Swarrow}{\rotatebox[origin=c]{45}{\(\Longleftarrow\)}}
\newcommand{\nc}{\newcommand}
\nc{\nt}{\newtheorem}
\nt{thm}{Theorem}[section]
\nt{cor}[thm]{Corollary}
\nt{prop}[thm]{Proposition}
\nt{lem}[thm]{Lemma}
\nt{assume}[thm]{Assumption}
\nt{claim}{Claim}

\theoremstyle{definition}
\nt{defn}[thm]{Definition}
\nt{rem}[thm]{Remark}
\nt{exa}[thm]{Example}

\nc{\ip}[2]{\mbox{$\langle #1,#2 \rangle$}}
\nc{\pf}{\noindent{\bf Proof\ \ }}
\nc{\finpf}{\hfill{$\Box$}\linespace}
\nc{\linespace}{\vspace{\baselineskip} \noindent}
\nc{\R}{{\mathbf R}}
\nc{\N}{{\mathbb N}}
\nc{\oR}{\overline{\R}}
\nc{\E}{{\bf E}}
\nc{\W}{{\mathcal W}}
\nc{\M}{\mathcal M}
\nc{\C}{\mathcal C}
\DeclareMathOperator{\CAT}{CAT}     
\DeclareMathOperator{\CBB}{CBB}     

\numberwithin{equation}{section}

\title{Local geometry of feasible regions via smooth paths}
\author{Adrian S. Lewis
\thanks{ORIE, Cornell University, Ithaca, NY. \texttt{people.orie.cornell.edu/aslewis} Research supported in part by National Science Foundation Grant DMS-2006990.
} 
\and 
Adriana Nicolae
\thanks{Department of Mathematics, Babe\c{s}-Bolyai University, 400084, Cluj-Napoca, Romania \newline
\texttt{anicolae@math.ubbcluj.ro}  Research supported in part by a grant of the Ministry of Research, Innovation and Digitization, CNCS/CCCDI -- UEFISCDI, project number PN-III-P1-1.1-TE-2019-1306, within PNCDI III.
} 
\and 
Tonghua Tian
\thanks{ORIE, Cornell University, Ithaca, NY.
\texttt{tt543@cornell.edu}}
}

\begin{document}
\date{}

\maketitle

\begin{abstract}
Variational analysis presents a unified theory encompassing in particular both smoothness and convexity.  In a Euclidean space, convex sets and smooth manifolds both have straightforward local geometry.  However, in the most basic hybrid case of feasible regions consisting of pre-images of convex sets under maps that are once (but not necessarily twice) continuously differentiable, the geometry is less transparent.  We define a new approximate convexity property, that holds both for such feasible regions and also for all prox-regular sets.  This new property requires that nearby points can always be joined by smooth feasible paths that are almost straight.  In particular, in the terminology of real algebraic geometry, such feasible regions are locally normally embedded in the Euclidean space.
\end{abstract}
\medskip

\medskip


\section{Introduction}
Variational analysis, as expounded in central texts such as \cite{VA,Mord_1,CLSW}, foregrounds a particular style of geometry, focused on normal and tangent vectors to sets $C$ in a Euclidean space  
$\R^n$.  The classical Bouligand tangent cone at a point  $x$  in  $C$,
\[
T_C(x) = \limsup_{t \downarrow 0} \frac{1}{t}(C-x),  
\]
plays a prominent role.  Variational analysis unifies two classical threads --- the study of convex sets  $D$ in a Euclidean space $\R^m$ on the one hand  and
${\mathcal C}^{(1)}$-smooth maps  $F \colon \R^n \to \R^m$ on the other --- leading to nonsmooth calculus rules for hybrid scenarios, central among them being the generalized equation
\[  
F(x) \in D.  
\]
Good behavior of this system at a solution $\bar x$  typically requires a certain ``constraint qualification'' involving the derivative map  $\nabla F(\bar x)$:
\begin{equation} \label{CQ}
T_D \big(F(\bar x)\big)  + \mbox{Range}\big(\nabla F(\bar x)\big) =   \R^m.
\end{equation}
Under this condition, in the terminology of \cite{VA}, the solution set  $C = F^{-1}(D)$ is called {\em amenable} at  $\bar x$, and {\em strongly amenable} if $F$ is in fact ${\mathcal C}^{(2)}$-smooth. Useful consequences of amenability include a formula for the tangent cone to  $C$  at  $\bar x$:
\[
T_{F^{-1}(D)}(\bar x) = \nabla F(\bar x)^{-1}\big(T_D(F(\bar x))\big).
\]  

Convex sets and ${\mathcal C}^{(1)}$-smooth manifolds, two basic building blocks of variational analysis, are everywhere amenable.  Convex sets arise by taking the map  $F$  to be the identity;  manifolds arise when the set  $D$  is a singleton.  When  $D$  is an orthant, the constraint qualification (\ref{CQ})  becomes the ``Mangasarian-Fromovitz'' condition central to classical nonlinear optimization.

Although the tools of variational analysis apply in principle to arbitrary closed sets  $C$  in  $\R^n$, some important avenues rely on regularity assumptions ensuring good behavior of all sets in some particular class.  Such assumptions are discussed extensively in \cite{VA}, amenability being one example.  From this perspective, rather than necessarily having in hand any particular representation   
$C = F^{-1}(D)$,  one simply assumes its existence and deduces useful variational consequences.  A second such assumption, common in various settings, is {\em prox-regularity} of  $C$  at a point  $\bar x$, which amounts to all points near $\bar x$  having unique nearest points in  $C$.  Amenability and prox-regularity are independent properties, neither implying the other, although both are implied by strong amenability.  In this work we propose a new candidate for a reasonable local assumption about sets in variational analysis, a property that we prove in our main result to be weaker than both amenability and prox-regularity.  We call this idea ``smooth approximate convexity'', a property we define shortly.

Among the regularity assumptions deployed throughout the variational analysis literature, the most basic is {\em Clarke regularity} \cite{Clarke75}, which for a closed set  $C$  at a point  $\bar x$  requires {\em inner semicontinuity} of the tangent cone: 
\[ 
T_C(\bar x) = \liminf_{x \to \bar x}  T_C(x).  
\]
Both amenability and prox-regularity imply Clarke regularity, and we prove that so does smooth approximate convexity.  To summarize, we will arrive at the following implications:
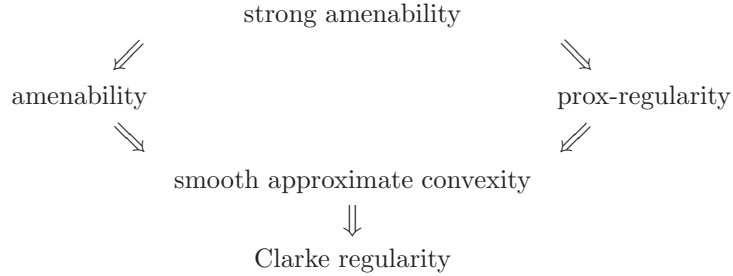
\begin{eqnarray*}
& \mbox{strong amenability} & \\ 
\Swarrow & & \Nwarrow \\ 
\mbox{amenability} & & 
\mbox{prox-regularity} \\ 
\Nwarrow & & \Swarrow \\ 
& \mbox{smooth approximate convexity} & \\ 
& \big\Downarrow & \\ 
& \mbox{Clarke regularity} &
\end{eqnarray*}

Basic variational analysis, despite its power, is surprisingly mute on the local metric geometry of sets.  Connectedness is an interesting example.  By definition, pairs of points in a convex set  $C$  are connected by line segments in  $C$.  Similarly, because any  ${\mathcal C}^{(1)}$-smooth manifold 
$C \subseteq \R^n$ is locally diffeomorphic to an open ball, nearby points are connected by  ${\mathcal C}^{(1)}$-smooth paths in $C$  that are close to line segments.  Both convex sets and ${\mathcal C}^{(1)}$-smooth manifolds are thus  {\em normally embedded} in the terminology of \cite{birbrair-mostowski}:  locally, the intrinsic distance between points in  $C$,  as measured by the length of rectifiable paths in  $C$  joining them, is bounded by a fixed multiple of their Euclidean distance.  This property is sometimes called ``Whitney 1-regularity'', due to its origins in \cite{whitney34}, or ``quasi-convexity'' in \cite{gromov} --- see also \cite{stasica,kurdyka92}.  Consider, then, the following natural question.\

\medskip
\noindent
{\bf Question 1.}  If a set in Euclidean space is amenable at a point, must it be normally embedded at that point?
\medskip

This basic question about the geometry of amenable sets seems unaddressed in the variational analysis literature.   The observation that amenability implies Clarke regularity does not help, because, for example, the set of points  $(u,v) \in \R^2$ satisfying $u^2 = v^4$  is Clarke regular at  $(0,0)$  but not normally embedded there.

The analogous question for prox-regular sets is easier to settle.  Prox-regularity implies that, locally,  $C$ has ``positive reach'' in the sense \cite{federer}.  By comparison with Clarke regularity, the metric geometry of sets of positive reach has a rich history, and indeed we shall see that prox-regular sets are normally embedded, and hence so are strongly amenable sets.

However, answering Question 1 by resorting to properties like prox-regularity or ${\mathcal C}^{(2)}$-smoothness seems unsatisfactory:   those properties lie beyond the question's two basic ingredients, namely convexity and ${\mathcal C}^{(1)}$-smoothness.  After all,  ${\mathcal C}^{(1)}$-smooth manifolds are normally embedded even though they may not be prox-regular: an example is the graph of the univariate function  $u \mapsto  |u|^{3/2}$,  which is not prox-regular at zero.  Both theory and computational practice \cite{bottou} encourage an approach relying on ${\mathcal C}^{(1)}$, rather than 
${\mathcal C}^{(2)}$-smoothness. 

We present an affirmative answer to Question 1: amenability implies normal embeddability. Indeed, we will not discuss normal embeddability further, because in fact we deduce a stronger consequence of amenability:  the smooth approximate convexity property, an idea that, like amenability itself, is framed simply in terms of ${\mathcal C}^{(1)}$-smoothness and a local convexity-like property.  However, unlike amenability, which is defined in terms of representability as an inverse image, smooth approximate convexity is an intrinsic geometric property.  

A set $C \subseteq \R^n$ is {\em smoothly approximately convex} at a point $\bar x \in C$ if, given any constant  $\epsilon > 0$, for all points  $x, x' \in C$  near  $\bar x$,  there exists a   
${\mathcal C}^{(1)}$-smooth map  $\eta$ from an open neighborhood of the interval $[0,1]$  into $\R^n$  satisfying  $\eta([0,1]) \subseteq C$, $\eta(0)= x$  and  $\eta(1)=x'$, and such that 
\[
\| \eta'(t) - (x'-x)\| \le \epsilon\|x'-x\| \qquad \mbox{for all}~ t \in [0,1].
\]
Geometrically, smooth approximate convexity of the set $C$ at the point $\bar x$ requires the existence, for each nearby pair of points $x$ and $x'$ in $C$, of a 
nearly straight ${\mathcal C}^{(1)}$-smooth path in $C$  connecting  $x$  and  $x'$.  Several more familiar properties follow immediately from smooth approximate convexity.  The restriction of the above path $\eta$ to $[0,1]$ has length no more than  $(1+\epsilon)\|x'-x\|$,  so the set  $C$  is normally embedded, as we claimed.  While we make no systematic effort to situate smooth approximate convexity precisely in the diffuse literature on general convexity-like properties, we do note that, along the path, the points  $\eta(t)$  for  $0 \le t \le 1$  satisfy the inequality
\[
\|\eta(t) - ((1-t)x + tx')\|  \le   2 \epsilon t(1-t)\|x-x'\|
\]
and so  $C$  is locally ``intrinsically approximately convex'' in the terminology of \cite{NgaiPenot07}.  Another immediate consequence is the ``uniform approximation by geodesics'' property introduced in \cite{lewis-lopez-nicolae-22} in the context of the method of alternating projections.

Our main aim in this work is speculative:  we propose smooth approximate convexity as a potential substitute for more familiar regularity conditions like amenability or prox-regularity.  We confine ourselves to one simple illustration from optimization.  Consider the problem of minimizing a  
${\mathcal C}^{(1)}$-smooth objective function  $f$  over a feasible region  $C \subseteq \R^n$.  The failure of a feasible point  $x$  to satisfy the standard first-order necessary condition amounts to the existence of a direction in the tangent cone  $T_C(x)$ whose inner product with the gradient  
$\nabla f(x)$  is negative.  Consequently, some sequence of feasible points  approaching  $x$  realizes a positive rate of objective decrease.  Computing such isolated points may be hard, unfortunately.  However, smooth approximate convexity guarantees instead a  ${\mathcal C}^{(1)}$-smooth feasible path from  $x$  along which the objective decreases linearly, a more promising computational scenario.  

Our discussion has focused on sets in Euclidean space, but much of our development extends naturally to subsets of Riemannian manifolds.

\section{Preliminaries}

\subsection{Euclidean spaces}

We fix some notation mostly following \cite{VA}. Let $\N$ denote the set of natural numbers excluding $0$. Given a Euclidean space $\R^n$, we denote by $\|\cdot\|$ the Euclidean norm on it. Denote the open ball and closed ball with center $x \in \R^n$ and radius $\delta > 0$ by $B_{\delta} (x)$ and $\overline{B}_{\delta} (x)$, respectively. A neighborhood of a point $x \in \R^n$ is an open set in $\R^n$ containing $x$. The range and kernel spaces of a linear mapping $A: \R^n \rightarrow \R^m$ are the sets $\mathrm{Range}(A) := \{A x : x \in \R^n\}$ and $\mathrm{Ker}(A) := \{x \in \R^n: Ax = 0\}$, respectively. Denote the graph and epigraph of a function $f : \R^n \rightarrow \R$ as $\mathrm{gph}\, f := \{(x, f(x)) : x \in \R^n\}$ and $\mathrm{epi}\, f := \{(x, t) : x \in \R^n, t \geq f(x)\}$.

Given a set $C \subset \R^n$ and points $x, y \in C$, a \textit{curve} in $C$ from $x$ to $y$ is a continuous mapping $\gamma: [a, b] \subseteq \R \rightarrow C$ such that $\gamma(a) = x$ and $\gamma(b) = y$. We say that $C$ is {\em connected} if for every pair of points $x, y \in C$ there is a curve in $C$ from $x$ to $y$. The set $C$ is {\em locally connected} at a point $\bar{x} \in C$ if for every open set $V$ containing $\bar x$, there exists an open set $U \subset V$ containing $\Bar{x}$ such that $C \cap U$ is connected. The set $C$ is called {\em locally closed} at $\bar{x}$ if $C \cap \overline{B}$ is closed for some closed ball $\overline{B}$ centered at $\bar{x}$. The {\em affine hull} of $C$, denoted by $\mathrm{aff}(C)$, is the smallest affine set in $\R^n$ containing $C$.

Let $\bar{x} \in C$. The {\em regular normal cone} to $C$ at $\Bar{x}$, denoted by $\widehat{N}_C (\Bar{x})$, is the set of all vectors $v$ such that
\begin{equation*}
    \langle v, x - \Bar{x} \rangle \leq o(\|x - \bar{x}\|) \qquad \text{for all}~ x \in C.
\end{equation*}
The \textit{normal cone} to $C$ at $\Bar{x}$, denoted by $N_C (\Bar{x})$, is the set of all vectors $v$ for which there exist two sequences $(x_k) \subseteq C$ and $(v_k) \subseteq \mathbb{R}^n$ with $v_k \in \widehat{N}_C (x_k)$ for all $k \in \mathbb{N}$ such that $x_k \rightarrow \Bar{x}$ and $v_k \rightarrow v$.

The notion of Clarke regularity discussed in the Introduction can equally be formulated using normal cones. The set $C$ is {\em Clarke regular} at $\Bar{x}$ if it is locally closed at $\Bar{x}$ and $N_C (\Bar{x}) = \widehat{N}_C (\Bar{x})$. The set $C$ is {\em prox-regular at $\Bar{x}$ for} $\bar{v} \in N_C (\Bar{x})$ if it is locally closed at $\Bar{x}$ and there exist constants $\epsilon, \rho > 0$ such that
\begin{equation*}
    \langle v, x' - x \rangle \leq \rho \|x' - x\|^2
\end{equation*}
for all $x, x' \in C \cap B_{\epsilon} (\Bar{x})$ and $v \in N_C (x) \cap B_{\epsilon} (\Bar{v})$; $C$ is {\em prox-regular} at $\Bar{x}$ if this property holds for all $\bar{v} \in N_C (\Bar{x})$. Prox-regularity is a stronger property than Clarke regularity. Prox-regularity of the set $C$ at $\Bar{x}$ can be equivalently characterized as the condition that every point around $\Bar{x}$ has a unique nearest point in $C$. This characterization will allow us to extend prox-regularity and related results to general metric spaces.

For $p \geq 0$, a $\C^{(p)}$-{\em smooth immersion} $H: U \rightarrow \R^n$, where $U$ is an open subset of $\R^k$, is a $\C^{(p)}$-smooth mapping such that $\nabla H(w)$ is injective for all $w \in U$. A $\C^{(p)}$-{\em smooth embedding} is a $\C^{(p)}$-smooth immersion that is also a homeomorphism onto its image. A $\C^{(p)}$-{\em smooth embedded submanifold} (in $\R^n$) is the image of a $\C^{(p)}$-smooth embedding. For a $\C^{(p)}$-smooth embedded submanifold $\mathcal{M}$ defined by the embedding $H: U \rightarrow \R^n$, its tangent space and normal space at a point $x \in \mathcal{M}$ are given by
\begin{equation*}
    T_{\mathcal{M}} (x) = \mathrm{Range}(\nabla H(w)) \quad \text{and} \quad N_{\mathcal{M}} (x) = \mathrm{Ker}(\nabla H(w)^{\top}),
\end{equation*}
where $w = H^{-1}(x)$. A mapping $F: \mathcal{M} \rightarrow \R$ is $\C^{(p)}$-smooth at $x$ if $F \circ H: U \rightarrow \R^m$ is $\C^{(p)}$-smooth at $w$ in the canonical sense, with its differential defined by $d F(x) = \nabla (F \circ H) (w) \nabla H(w)^{-1} : T_{\mathcal{M}} (x) \rightarrow \R^m$. If $F$ extends to a $\C^{(p)}$-smooth mapping $\widetilde{F}$ defined on a neighborhood of $x$ in $\R^n$, then $d F(x) = \nabla \widetilde{F} (x) |_{T_{\mathcal{M}} (x)}$.

\subsection{Metric spaces}

In this paper, we also work with more general metric spaces $(X, d)$. (When $X$ is a Euclidean space, $d(x, y) = \|x - y\|$.) At this point we only give a few basic definitions, but we will include a more detailed discussion in Appendix \ref{sec:equi}.

Still denote the open ball and closed ball with center $x \in X$ and radius $\delta > 0$ by $B_{\delta} (x)$ and $\overline{B}_{\delta} (x)$, respectively. The distance from a point $x \in X$ to a set $C \subseteq X$ is $d(x, C) := \inf_{y \in C} d(x, y)$, and the projection of $x$ onto $C$ is the possibly empty set $P_C (x) := \{y \in X : d(x, y) = d(x, C)\}$.

Given a {\em curve} in $X$, namely a continuous mapping $\gamma: [a, b] \subset \R \rightarrow X$, define its {\em length} as
\begin{equation*}
    L(\gamma) := \sup\big\{\sum_{i=0}^{n-1} d(\gamma (t_i), \gamma(t_{i+1})) : a=t_0\leq t_1\leq \cdots \leq t_n=b, n \in \N\big\}.
\end{equation*}
For every set $C \subseteq X$ and points $x, y \in C$, define the {\em intrinsic distance} $d_C (x, y) \in [0, \infty]$ in $C$ between $x$ and $y$ as the infimum of the lengths of the curves in $C$ from $x$ to $y$. We say $(X, d)$ is a {\em length space} if $d_X = d$.

Let $x, y \in X$. A {\em geodesic} from $x$ to $y$ is a curve $\gamma : [0, L(\gamma)] \rightarrow X$ from $x$ to $y$ such that
\begin{equation*}
    d\big(\gamma(t_1), \gamma(t_2)\big) = |t_1 - t_2| \qquad \text{for all}~ t_1, t_2 \in [0, L(\gamma)].
\end{equation*}
The image $\gamma([0,L(\gamma)])$ is called a {\em geodesic segment}. We say that $X$ is a {\em (uniquely) geodesic space} if for every pair of points $x, y \in X$ there exists a (unique) geodesic from $x$ to $y$. Every geodesic space is a length space. Conversely, by the Hopf-Rinow Theorem, every length space that is complete and locally compact is geodesic.

In this paper, it is often easier to work with reparametrized geodesics as follows. Given $x, y \in X$, an {\em averaging map} from $x$ to $y$ is a curve $\gamma : [0, 1] \rightarrow X$ from $x$ to $y$ such that
\begin{equation}\label{eq:def-averaging-map}
    d\big(\gamma(t_1), \gamma(t_2)\big) = |t_1 - t_2| \cdot d(x, y) \qquad \text{for all}~ t_1, t_2 \in [0, 1].
\end{equation}
Obviously, a curve $\gamma : [0, 1] \rightarrow X$ is an averaging map if and only if its reparametrization $t \mapsto \gamma \big( t/d(x, y) \big)$ is a geodesic.

\section{Smooth approximate convexity in Euclidean spaces}
\label{sec:eclidean}

Approximate convexity has been discussed extensively in the literature. In particular, \cite{NgaiPenot07} calls a set $C \subseteq \R^n$ ``intrinsically approximately convex'' at $\bar{x} \in C$ if given any $\epsilon > 0$, all $x, x' \in C$ around $\bar{x}$ satisfy
\begin{equation}\label{eq:def-intrinsic}
    d ((1-t)x + tx', C) \leq \epsilon t (1-t) \|x'-x\| \qquad \mbox{for all}~ t \in [0,1].
\end{equation}
Analogously, \cite{NLT2000} defines a notion of approximate convexity for functions. A function $f: \R^n \rightarrow \R$ is called ``approximately convex'' at $\Bar{x}$ if given any $\epsilon > 0$, all $x, x'$ around $\Bar{x}$ satisfy
\begin{equation}\label{eq:def-fac}
    f((1-t)x + tx') \leq (1-t) f(x) + t f(x') + \epsilon t(1-t) \|x'-x\| \qquad \mbox{for all}~ t \in [0,1].
\end{equation}

 Based on these ideas, we introduce a new definition that not only captures approximate convexity but in a $\C^{(1)}$-smooth sense. Given a set $C \subseteq \R^n$, as in \cite{Lee}, we call a curve $\gamma: [a, b] \rightarrow C$ a \textit{smooth curve segment in $C$}  if it has an extension to a $\C^{(1)}$-smooth map from a neighborhood of $[a, b]$ to $\R^n$. (Note that the mapping $\eta$ used in the Introduction to define smooth approximate convexity is such an extension.)

\begin{defn}\label{def:sac}
    A set $C \subseteq \R^n$ is {\em smoothly approximately convex} at a point $\Bar{x} \in C$ if given any $\epsilon > 0$, for all $x, x' \in C$ around $\Bar{x}$, there is a smooth curve segment $\gamma: [0, 1] \rightarrow C$ with $\gamma(0) = x$ and $\gamma(1) = x'$ that satisfies
    \begin{equation}\label{eq:def-sac}
        \|\gamma' (t) - (x' - x) \| \leq \epsilon \|x'-x\| \qquad \mbox{for all}~ t \in [0, 1].
    \end{equation}
    We call such a map $\gamma$ an {\em $\epsilon$-path} from $x$ to $x'$. Moreover, we say $C$ is smoothly approximately convex if it is smoothly approximately convex at every $\Bar{x} \in C$.
\end{defn}

\begin{rem}\label{rmk:eps-path}
We can immediately see that the length of any $\epsilon$-path $\gamma$ from $x$ to $x'$ is at most $(1+\epsilon)\|x'-x\|$. Moreover, 
	\begin{equation}\label{eq:rmk-eps-path}
		\|\gamma(t) - \left((1-t) x + t x'\right)\|  = \left\| \int_0^t \left(\gamma' (\tau) - (x'-x)\right) \mathrm{d}\tau \right\| \le 		\epsilon t \|x'-x\|  \qquad \mbox{for all}~ t \in [0, 1].
	\end{equation}
By symmetry of $\gamma$ with respect to $x$ and $x'$, we can deduce that 
\[\|\gamma(t) - \left((1-t) x + t x'\right)\|  \le 2\epsilon t(1-t) \|x'-x\|  \qquad \mbox{for all}~ t \in [0, 1].\]

\end{rem}

\begin{rem}\label{rmk:fac}
    If a function $f$ is approximately convex at $\bar{x}$, then given any $\bar{s} \geq f(\bar{x})$, by considering, for $(x, s), (x', s') \in \mathrm{epi}\, f$ around $(\bar{x}, \bar{s})$, the curve
    \begin{equation*}
        \gamma(t) = \left((1-t)x + t x', (1-t) s + t s' + \epsilon t(1-t) \|x'-x\|\right)\qquad \mbox{for all}~ t \in [0,1],
    \end{equation*}
   we can see that $\mathrm{epi}\, f$ is smoothly approximately convex at $(\bar{x}, \bar{s})$. Moreover, it is shown in \cite{DG04} that a locally Lipschitz function is lower-$\C^{(1)}$ if and only if it is approximately convex. Hence, the epigraph of any lower-$\C^{(1)}$ function is smoothly approximately convex. On the other hand, it is easy to see that smooth approximate convexity implies intrinsic approximate convexity.
\end{rem}

Henceforth, we will concentrate on the idea of smooth approximate convexity, and will not discuss further intrinsic approximate convexity or approximate convexity. Clearly, the graph and epigraph of any $\C^{(1)}$-smooth function are smoothly approximately convex. Another family of functions whose epigraphs are smoothly approximately convex is described in the following proposition.

\begin{prop}[Smooth approximate convexity of epigraphs]
\label{prop:0lip}
    If a function $f: \R^n \rightarrow \R$ is locally Lipschitz continuous at $\bar{x}$ with local Lipschitz modulus $0$, then $\mathrm{epi}\, f$ is smoothly approximately convex at $(\bar{x}, \bar{s})$ for all $\bar{s} \geq f(\bar{x})$.
\end{prop}

\pf
    Fix any $\epsilon > 0$. There exists a neighborhood $U$ of $\bar{x}$ such that $f$ is $\frac{\epsilon}{2}$-Lipschitz continuous on $U$. Then for all $x, x' \in U$ we have
    \begin{align}
        f((1-t)x + tx') &\leq f(x) + \frac{\epsilon}{2} \|t(x' - x)\|, \label{eq:0lip-1} \\
         f((1-t)x + tx') &\leq f(x') + \frac{\epsilon}{2} \|(1-t)(x' - x)\|. \label{eq:0lip-2}
    \end{align}
    Multiplying (\ref{eq:0lip-1}) with $1-t$ and (\ref{eq:0lip-2}) with $t$ and adding up, we obtain (\ref{eq:def-fac}), which means $f$ is approximately convex at $\bar{x}$. The conclusion follows from Remark \ref{rmk:fac}.
\finpf

The following lemma reveals that $\C^{(1)}$-smooth embeddings preserve smooth approximate convexity, which is essential for our proofs later.

\begin{lem}[Smooth embeddings preserve smooth approximate convexity]
\label{lem:sac-img}
    Suppose the set $C \subseteq \R^n$ is smoothly approximately convex at $\Bar{x} \in C$. If the map $H: U \rightarrow \R^m$, where $U$ is a neighborhood of $\Bar{x}$ in $\R^n$, is a $\C^{(1)}$-smooth embedding, then $H(C \cap U)$ is smoothly approximately convex at $H(\Bar{x})$.
\end{lem}

\pf
    Let $\Bar{y} = H(\Bar{x})$ and $\mathcal M = H(U)$. Suppose $H = (h_1, \dots, h_m)$. Fix any $\epsilon' > 0$. By the injectivity of $\nabla H(\Bar{x})$ and $\C^{(1)}$-smoothness of $H$, we find a neighborhood $\widetilde{U} \subseteq U$ of $\Bar{x}$ in $\R^n$ and constants $c_1, c_2 > 0$ such that for all $x, x' \in \widetilde{U}$,
  \[ \|H(x) - H(x')\| \geq c_1 \|x - x'\|\]
and       
    \[ \|\nabla h_i(x)\| \leq \frac{c_2}{m}, \qquad \|\nabla h_i(x) - \nabla h_i(x')\| \leq \frac{c_1 \epsilon'}{2 m} \qquad \mbox{for all}~ i \in\{1, \dots, m\}.\]

Let $\epsilon = \min\{\frac{c_1}{2 c_2} \epsilon', \frac{1}{2}\}$. Since $C$ is smoothly approximately convex at $\Bar{x}$, there exists $\delta > 0$ such that $B_{2\delta} (\Bar{x}) \subseteq \widetilde{U}$ and for all $x, x' \in B_{\delta}(\Bar{x}) \cap C$, there is an $\epsilon$-path from $x$ to $x'$. Since $H: U \rightarrow \mathcal M$ is a homeomorphism, there exists $\delta' > 0$ such that
    \begin{equation*}
        H^{-1} (B_{\delta'}(\Bar{y}) \cap \mathcal M) \subseteq B_{\delta}(\Bar{x}).
    \end{equation*}
    
    Fix any $y, y' \in B_{\delta'}(\Bar{y}) \cap H(C \cap U)$. Then there exist $x, x' \in B_{\delta}(\Bar{x})$ such that $y = H(x)$ and $y' = H(x')$. Moreover, the injectivity of $H$ implies $x, x' \in C$. Let $\gamma$ be an $\epsilon$-path from $x$ to $x'$. Using \eqref{eq:rmk-eps-path}, we deduce that
    \[
    \|\gamma(t) - \Bar{x}\| \le \|\gamma(t) - \left((1-t) x + t x'\right)\| + \|\left((1-t) x + t x'\right) - \Bar{x}\| \le \epsilon t \|x'-x\| + \delta < 2\epsilon \delta + \delta \le 2\delta,
    \]
    for all $t \in [0,1]$. Consequently, $\gamma ([0,1]) \subseteq B_{2\delta}(\Bar{x}) \subseteq \widetilde{U}$. Let
    \begin{equation*}
        \Tilde{\gamma} = H \circ \gamma: [0, 1] \rightarrow H(C\cap U).
    \end{equation*}
    For each $i$, by the mean value theorem, there exists $\xi_i \in \widetilde{U}$ such that
    \begin{equation*}
        h_i(x') - h_i(x) = \nabla h_i(\xi_i)^{\top} (x' - x).
    \end{equation*}
    Hence,
    \begin{align*}
        \|\Tilde{\gamma}'(t) - (y'-y)\|
        &\leq \sum_{i = 1}^m \left\|\nabla h_i (\gamma(t))^\top \gamma'(t) - \nabla h_i(\xi_i)^{\top} (x' - x) \right\| \\
        &\leq m \left(\frac{c_2}{m} \epsilon \|x'-x\| + \frac{c_1 \epsilon'}{2 m} \|x'-x\|\right) \\
        &\leq \left(\frac{c_2 \epsilon}{c_1} + \frac{\epsilon'}{2} \right) \|y' - y\| \leq \epsilon' \|y' - y\| \qquad \mbox{for all}~ t \in [0, 1].
    \end{align*}
\finpf

\subsection{Other regularity conditions}

Various regularity conditions have been introduced in the literature to ensure good behavior of sets in different aspects. Conditions that appear unrelated sometimes have surprising implication relationships. It is worthwhile to inspect the relationships between smooth approximate convexity and common regularity conditions, including Clarke regularity and prox-regularity. As we will see, for closed sets, smooth approximate convexity is strictly stronger than Clarke regularity and strictly weaker than prox-regularity.

The proof for prox-regularity will be presented in Section \ref{sec:ecli-prox}. For now, we demonstrate the relationship between smooth approximate convexity and Clarke regularity through two other conditions, which are also interesting for their own sake. The first one, super-regularity, was introduced by Lewis, Luke, and Malick in \cite{lewis-luke-malick} as a condition for local linear convergence of alternating projections. A set $C \subseteq \R^n$ is ``super-regular'' at a point $\bar{x} \in C$ if it is locally closed there and given any $\epsilon > 0$, all $x, x' \in C$ around $\bar{x}$ satisfy
\begin{equation}\label{eq:def-super}
    \langle v, x' - x \rangle \leq \epsilon \|v\| \cdot \|x'-x\|\qquad \mbox{for all}~ v \in N_C (x).
\end{equation}
Later, Lewis, L\'opez-Acedo and Nicolae \cite{lewis-lopez-nicolae-22} extended the analysis of the local linear convergence of alternating projections to the setting of metric spaces with bounded curvature by studying the second condition, the notion of uniform approximation by geodesics. When restricting the setting to $\R^n$, a set $C$ is ``uniformly approximable by geodesics'' (UAG) at a point $\bar{x} \in C$ if given any $\epsilon > 0$, for all distinct $x, x' \in C$ around $\bar{x}$, there is a map $\gamma: [0, 1] \rightarrow C$ and a direction $d \in \R^n \backslash \{0\}$ such that $\gamma(0) = x, \gamma(1) = x'$ and
\begin{equation}\label{eq:def-uag}
    \|\gamma (t) - (x + td)\| \leq \epsilon t \|d\| \qquad \mbox{for all}~ t \in [0, 1].
\end{equation}
The authors showed in \cite{lewis-lopez-nicolae-22} that if a closed set is UAG at a point, then it is super-regular there. Furthermore, using \eqref{eq:rmk-eps-path}, we can easily see that smooth approximate convexity implies the UAG property.

\begin{prop}[Implications of smooth approximate convexity]
\label{prop:implications}
    If a set $C \subseteq \R^n$ is smoothly approximately convex at $\Bar{x} \in C$, then it is UAG there. If it is furthermore locally closed at $\Bar{x}$, then it is super-regular and Clarke regular there.
\end{prop}

Next we present two examples. The first example illustrates that an everywhere Clarke regular set might not be super-regular, and therefore neither smoothly approximately convex. The second example from \cite{lewis-lopez-nicolae-22} shows that the UAG property at a point does not ensure smooth approximate convexity there, therefore neither does super-regularity.

\begin{exa}[Clarke regularity does not imply smooth approximate convexity]
    Consider the closed set
    \begin{equation*}
        C = \{(x, y) \in \R^2 : y = x^2 \text{ or } y = - x^2\}.
    \end{equation*}
    It is easy to see that
    \[N_C (x, x^2) = \widehat{N}_C (x, x^2) = (-2x, 1) \cdot \R \quad \text{and} \quad N_C (x, -x^2) = \widehat{N}_C (x, -x^2) = (2x, 1) \cdot \R,\]
    for all $x \in \R$. So $C$ is everywhere Clarke regular. To see that it is not super-regular at $(0, 0)$, consider the points $x_k = (\frac{1}{k}, \frac{1}{k^2}), x_k' = (\frac{1}{k}, -\frac{1}{k^2})$ and $v_k = (\frac{2}{k}, -1) \in N_C (x_k)$ for all $k \geq 1$. Clearly, $x_k, x_k' \rightarrow (0, 0)$ as $k \rightarrow \infty$. However, for each $k$ we have
    \begin{equation*}
        \langle v_k, x_k' - x_k \rangle = \frac{2}{k^2} > \frac{1}{3} \|v_k\| \cdot \|x_k' - x_k\|.
    \end{equation*}
    Hence (\ref{eq:def-super}) fails for $\epsilon = \frac{1}{3}$. By Proposition \ref{prop:implications}, $C$ is not smoothly approximately convex at $(0, 0)$.
\end{exa}

\begin{exa}[UAG does not imply smooth approximate convexity]
\label{ex:super-sac}
    Define the function $f: \R \rightarrow \R$ as follows
    \begin{equation*}
        f(x) = 
        \begin{cases} \displaystyle
            \frac{1}{2^k} \Big( x - \frac{1}{2^{k+1}}\Big) & \displaystyle \text{ if } x \in \Big(\frac{1}{2^{k+1}}, \frac{3}{2^{k+2}}\Big], k \in \N \\[8pt]  \displaystyle
            \frac{1}{2^k} \Big(\frac{1}{2^{k}} - x\Big) & \displaystyle \text{ if } x \in \Big(\frac{3}{2^{k+2}}, \frac{1}{2^{k}}\Big], k \in \N \\[8pt] 
            0 & \text{otherwise.}
        \end{cases}
    \end{equation*}
    Let $C = \mathrm{gph}\, f$. It was proved in \cite{lewis-lopez-nicolae-22} that $C$ is UAG at $(0, 0)$. For the sake of contradiction, assume it is smoothly approximately convex there. Then for $\epsilon = \frac{1}{2}$ and a sufficiently large $k$, there is an $\epsilon$-path $\gamma = (\gamma_1, \gamma_2)$ from $x = (0, 0)$ to $x' = (\frac{1}{2^k}, 0)$. By (\ref{eq:def-sac}) we can deduce that
    \begin{equation*}
        |\gamma_1'(0) - \frac{1}{2^k}| \leq \|\gamma'(0) - (x'-x)\| \leq \epsilon \|x'-x\| < \frac{1}{2^k},
    \end{equation*}
    which means $\gamma_1'(0) > 0$. So there is a $\C^{(1)}$-smooth map $g: (0, \delta) \rightarrow (0,\infty)$ for some $\delta > 0$ such that $\gamma_1 \circ g = \mathrm{id}_{(0, \delta)}$. Hence we have
    \begin{equation*}
        \gamma_2 \circ g = f \circ \gamma_1 \circ g = f|_{(0, \delta)}.
    \end{equation*}
    It follows that $f$ is $\C^{(1)}$-smooth on $(0, \delta)$, which is obviously not true.
\end{exa}

\begin{rem}\label{rmk:sac-nonlocal}
    In Example \ref{ex:super-sac}, it is easy to see that $f$ is locally Lipschitz continuous at $0$ with local Lipschitz modulus $0$. So by Proposition \ref{prop:0lip}, we know that the set $\mathrm{epi}\, f$ is smoothly approximately convex at $(0, 0)$. However, it is obviously not Clarke regular at $(\frac{3}{2^{k+2}}, \frac{1}{2^{2k+2}})$ for all $k \in \N$, therefore not smoothly approximately convex at these points either. This means smooth approximate convexity is not persistent nearby.
\end{rem}

We conclude this section by mentioning another regularity property, called subsmoothness, which can be regarded as a translation of prox-regularity in a first-order context. Namely, a set $C \subseteq \R^n$ is ``subsmooth'' at a point $\bar{x} \in C$ if it is locally closed there and given any $\rho > 0$, there exists $\epsilon > 0$ such that $\langle v, x' - x \rangle \leq \rho \|x' - x\|$ for all $x, x' \in C \cap B_{\epsilon} (\Bar{x})$ and all $v \in N^c_C (x) \cap \overline{B}_{1} (0)$, where $N^c_C (x)$ is the closed convex hull of $N_C (x)$ called the Clarke normal cone to $C$ at $x$. Subsmoothness is strictly weaker than prox-regularity and strictly stronger than super-regularity, intrinsic approximate convexity, and Clarke regularity. Moreover, a function $f : \R^n \to \R$ that is locally Lipschitz continuous at $\bar{x} \in \R^n$ is approximately convex at $\bar{x}$ if and only if its epigraph is subsmooth at $(\bar{x},f(\bar{x}))$. We refer to the survey \cite{thibault} for a comprehensive study of subsmooth sets.

\subsection{Amenability vs. smooth approximate convexity}

The first main result in this paper states that all amenable sets are smoothly approximately convex. To see this, we start with some discussions about the structure of amenable sets. Note that one can equally reformulate the constraint qualification \eqref{CQ} defining amenable sets as follows.

\begin{defn}\label{def:amenable}
    A set $C \subseteq \R^n$ is {\em amenable} at $\Bar{x} \in C$ if there is a neighborhood $V$ of $\Bar{x}$ along with a $\C^{(1)}$ mapping $F$ from $V$ into a space $\R^m$ and a closed, convex set $D\subseteq \R^m$ such that
    \begin{equation*}
        C \cap V = \left\{x \in V : F(x) \in D \right\}
    \end{equation*}
    and the following constraint qualification is satisfied:
    \begin{equation}\label{eq:CQ}
        \text{the only vector}~ y \in N_D (F(\Bar{x})) ~\text{with}~ \nabla F(\Bar{x})^{\top}y = 0 ~\text{is}~ y=0.
    \end{equation}
    In this case, we call $(V, F, D)$ a {\em representation} of the set $C$ around $\Bar{x}$.
\end{defn}

To prove the desired implication, we start with a special type of amenable sets for which the convex sets in the representation have nonempty interiors.

\begin{lem}\label{lem:CQ}
    In Definition \ref{def:amenable}, if $\mathrm{int}(D) \neq \emptyset$, then the constraint qualification (\ref{eq:CQ}) is equivalent to the following condition:
    \begin{equation}\label{eq:CQ-dual}
       \text{there exists}~ r \in \R^n ~\text{such that}~ F(\Bar{x}) + \nabla F(\Bar{x}) r \in \mathrm{int}(D).
    \end{equation}
\end{lem}

\pf
    The fact that $(\ref{eq:CQ-dual})$ implies $(\ref{eq:CQ})$ follows since $\langle y, z - F(\Bar{x}) \rangle < 0$ for all $y \in N_D (F(\Bar{x})) \backslash \{0\}$ and $z \in \mathrm{int}(D)$. To see the converse, assume (\ref{eq:CQ-dual}) is not true. Then we have
    \begin{equation*}
        0 \notin \mathrm{Range}(\nabla F(\Bar{x})) - \left(\mathrm{int}(D) - F(\Bar{x})\right).
    \end{equation*}
    Applying \cite[Theorem 2.39]{VA}, $\mathrm{Range}(\nabla F(\Bar{x}))$ and $\mathrm{int}(D) - F(\Bar{x})$ can be separated by a hyperplane, namely there exist $a \in \R^m \backslash \{0\}$ and $\alpha \in \R$ such that
    \begin{align*}
        \langle a, v\rangle &\geq \alpha\qquad \mbox{for all}~ v \in \mathrm{Range}(\nabla F(\Bar{x})), \\
        \langle a, v\rangle &\leq \alpha\qquad \mbox{for all}~ v \in \mathrm{int}(D) - F(\Bar{x}).
    \end{align*}
    Clearly, it must hold that $\alpha = 0$, $a \in \mathrm{Ker}(\nabla F(\Bar{x})^{\top})$ and
    \begin{equation*}
        \langle a, v - F(\Bar{x})\rangle \leq 0\qquad \mbox{for all}~ v \in D,
    \end{equation*}
    which means $a \in N_D (F(\Bar{x}))$. This contradicts (\ref{eq:CQ}).
\finpf

\begin{prop}\label{prop:amen-int}
    Let $C \subseteq \R^n$ be an amenable set with representation $(V, F, D)$ around $\Bar{x} \in C$. Suppose $\mathrm{int}(D) \neq \emptyset$. Then $C$ is smoothly approximately convex at $\Bar{x}$. Moreover, the $\epsilon$-path $\gamma$ connecting any two points in $C$ around $\Bar{x}$ can be taken so that $\gamma((0,1)) \subseteq \mathrm{int}(C)$.
\end{prop}

\pf
    If $F(\Bar{x}) \in \mathrm{int}(D)$, then $\Bar{x} \in \mathrm{int}(C)$, therefore $C$ is smoothly approximately convex at $\Bar{x}$ trivially. Henceforth we assume $F(\Bar{x}) \notin \mathrm{int}(D)$. By Lemma \ref{lem:CQ}, there exists $w \in \R^n$ and $\lambda > 0$ such that $\|w\|=1$ and
    \begin{equation*}
        F(\Bar{x}) + \lambda \nabla F(\Bar{x}) w \in \mathrm{int}(D).
    \end{equation*}
    So there exists $\delta > 0$ such that
    \begin{equation*}
        B_{\delta} \left(F(\Bar{x}) + \lambda \nabla F(\Bar{x}) w\right) \subseteq \mathrm{int}(D).
    \end{equation*}
    Fix any $\epsilon \in (0, 1]$. Suppose $F(x) = (f_1(x), \dots, f_m(x))$ for all $x$ close to $\bar x$. By the $\C^{(1)}$-smoothness of $F$, there exists $\delta' \in (0, \frac{4 \lambda}{\epsilon})$ such that for all $x \in B_{\delta'}(\Bar{x})$, $F$ is defined at $x$ with
    \begin{equation*}
        \|F(x) - F(\Bar{x})\| < \frac{\delta}{2} \quad\text{and}\quad \|\nabla f_i (x) - \nabla f_i (\Bar{x})\| < \frac{\epsilon \delta}{16 m \lambda}\qquad \mbox{for all}~ i\in\{1, \dots, m\}.
    \end{equation*}
    Fix any distinct $x, x' \in B_{\delta'/2} (\Bar{x}) \cap C$. Consider the map $\gamma: [0, 1] \rightarrow \R^n$ defined as
    \begin{equation*}
        \gamma (t) = (1-t)x + tx' + \epsilon t (1-t) \|x - x'\|w, \quad t \in [0, 1].
    \end{equation*}
    Obviously, $\gamma$ satisfies (\ref{eq:def-sac}). We only need to show that $\gamma(t) \in \mathrm{int}(C)$ for all $t \in (0, 1)$.
    
    Fix any $t \in (0, \frac{1}{2}]$. By the mean value theorem, for each $i \in \{1, \ldots, m\}$ there exist $\xi_i \in [x, \gamma(t)]$ and $\eta_i \in [x, x']$ such that
    \begin{align*}
        f_i(\gamma (t)) &= f_i(x) + \nabla f_i(\Bar{x})^{\top} \left(tx' - tx + \epsilon t (1-t) \|x - x'\|w\right) + h_i(t), \\
        f_i(x') &= f_i(x) + \nabla f_i(\Bar{x})^{\top} \left(x' - x\right) + g_i(t),
    \end{align*}
    where
    \begin{align*}
        h_i(t) &= (\nabla f_i(\xi_i) - \nabla f_i(\Bar{x}))^{\top} \left(tx' - tx + \epsilon t (1-t) \|x - x'\|w\right), \\
        g_i(t) &= (\nabla f_i(\eta_i) - \nabla f_i(\Bar{x}))^{\top} \left(x' - x\right).
    \end{align*}
    Note that we have $x, x' \in B_{\delta'/2} (\Bar{x})$ and
    \begin{equation*}
        \|\gamma (t)-\Bar{x}\| \leq (1-t)\|x-\Bar{x}\| + t\|x'-\Bar{x}\| + \frac{\epsilon}{4} \|x - x'\| < \delta',
    \end{equation*}
    therefore $\xi_i, \eta_i \in B_{\delta'} (\Bar{x})$. Hence,
    \begin{align*}
        |h_i(t)| &\leq \left\|\nabla f_i(\xi_i) - \nabla f_i(\Bar{x})\right\| \cdot \left\|tx' - tx + \epsilon t (1-t) \|x - x'\|w\right\| \\
        &\leq \frac{\epsilon \delta}{16 m \lambda} \cdot 2t \|x'-x\| = \frac{\epsilon \delta}{8 m \lambda} t \|x'-x\|, \\
        |g_i(t)| &\leq \left\|\nabla f_i(\eta_i) - \nabla f_i(\Bar{x})\right\| \cdot \left\|x' - x\right\| \leq \frac{\epsilon \delta}{16 m \lambda} \|x'-x\|.
    \end{align*}
    Let $h(t) = (h_1(t), \dots, h_m(t))$ and $g(t) = (g_1(t), \dots, g_m(t))$. Then
    \begin{align*}
        F(\gamma(t)) &= F(x) + \nabla F(\Bar{x}) \left(tx' - tx + \epsilon t (1-t) \|x - x'\|w\right) + h(t), \\
        F(x') &= F(x) + \nabla F(\Bar{x}) \left(x' - x\right) + g(t),
    \end{align*}
    where
    \[\|h(t)\| \leq \frac{\epsilon \delta}{8 \lambda} t \|x'-x\| \quad \text{and} \quad \|g(t)\| \leq \frac{\epsilon \delta}{16 \lambda} \|x'-x\|.
    \end{equation*}
    Hence,
    \begin{equation*}
        F(\gamma (t)) = (1-t) F(x) + t F(x') + \epsilon t (1-t) \|x - x'\| \nabla F(\Bar{x})w + h(t) - t g(t),
    \end{equation*}
    where
    \begin{equation*}
        \|h(t) - t g(t)\| < \frac{\epsilon \delta}{4 \lambda} t \|x'-x\| \leq \frac{\epsilon \delta}{2 \lambda} t (1-t) \|x'-x\|.
    \end{equation*}
    Letting $\mu = \frac{\epsilon}{\lambda} t (1-t) \|x-x'\|$, we have $\mu \in (0, 1)$ and
    \begin{equation*}
        \epsilon t (1-t) \|x - x'\| \nabla F(\Bar{x})w + h(t) - t g(t) \in \mu \cdot B_{\delta/2} (\lambda \nabla F(\Bar{x})w).
    \end{equation*}
    Let $y = (1-t) F(x) + t F(x')$. Then $y \in B_{\delta/2} (F(\Bar{x})) \cap D$. Thus,
    \begin{align*}
        F(\gamma(t)) &\in y + \mu \cdot B_{\delta/2} (\lambda \nabla F(\Bar{x})w) \\
        &= \left(1 - \mu\right) y + \mu \cdot \left(y + B_{\delta/2} (\lambda \nabla F(\Bar{x})w)\right) \\
        &\subseteq \left(1 - \mu\right) y + \mu \cdot \left(B_{\delta/2} (F(\Bar{x})) + B_{\delta/2} (\lambda \nabla F(\Bar{x})w)\right) \\
        &\subseteq \left(1 - \mu\right) y + \mu \cdot B_{\delta} \left(F(\Bar{x}) + \lambda \nabla F(\Bar{x})w\right) \\
        &\subseteq \left(1 - \mu\right) \cdot D + \mu \cdot \mathrm{int}(D) \subseteq \mathrm{int}(D).
    \end{align*}
    So $\gamma$ maps $(0, \frac{1}{2}]$ into $\mathrm{int}(C)$. By symmetry, it follows that $\gamma$ also maps $[\frac{1}{2}, 1)$ into $\mathrm{int}(C)$.
\finpf

\begin{rem}\label{rmk:connected}
    We can see from Proposition \ref{prop:amen-int} that $F^{-1} (\mathrm{int}(D))$ is a connected set around $\Bar{x}$.
\end{rem}

Next we use Lemma \ref{lem:sac-img} to extend the above result to general amenable sets.

\begin{thm}[Amenability implies smooth approximate convexity]
\label{thm:euclidean}
    If a set $C \subseteq \R^n$ is amenable at $\Bar{x} \in C$, then it is smoothly approximately convex at $\Bar{x}$.
\end{thm}

\pf
    Let $(V, F, D)$ be a representation of $C$. By Proposition \ref{prop:amen-int}, we only need to consider the case when $\mathrm{int}(D) = \emptyset$. Let $K = \mathrm{aff}(D)$. We can assume $K = \mathrm{span}\{e_1, \dots, e_k\}$, where $k < m$ and $e_i$ denotes the $i$th standard basis vector in $\R^n$. Suppose $F(x) = (f_1(x), \dots, f_m(x))$, $x \in V$. Define the mappings $F_1$ and $F_2$ as
    \begin{equation*}
        F_1(x) = (f_1 (x), \dots, f_k (x)), ~~ F_2(x) = (f_{k+1} (x), \dots, f_m (x)), \quad x \in V.
    \end{equation*}
    Then $F^{-1}(K) = F_2^{-1}(0)$. Note that for any $u \in \mathrm{Ker}(\nabla F_2(\Bar{x})^{\top})$, we have
    \begin{equation*}
        (0, u) \in N_D(F(\Bar{x})) \cap \mathrm{Ker}(\nabla F(\Bar{x})^{\top}) = \{0\},
    \end{equation*}
    therefore $\nabla F_2(\Bar{x})$ is surjective. Thus, $F_2^{-1}(0)$ is a $\C^{(1)}$ manifold around $\Bar{x}$ of dimension $n-m+k$. Hence there exists a neighborhood $U$ of $0$ in $\R^{n-m+k}$, a neighborhood $\widetilde{V} \subseteq V$ of $\Bar{x}$ in $\R^n$ and a $\C^{(1)}$-smooth embedding $H: U \rightarrow \R^n$ with $H(0) = \Bar{x}$ such that
    \begin{equation*}
        F^{-1}(K) \cap \widetilde{V} = F_2^{-1}(0) \cap \widetilde{V} = H(U).
    \end{equation*}
    Define the map $P: \R^m \rightarrow \R^k$ as
    \begin{equation*}
        P(y_1, \dots, y_m) = (y_1, \dots, y_k).
    \end{equation*}
    Let $\widetilde{F} = F_1 \circ H$ and $\widetilde{D} = P(D)$. Then we have
    \begin{align}
        C \cap \widetilde{V} &= \{x \in \widetilde{V}: F(x) \in K \text{ and } F_1(x) \in \widetilde{D}\} \nonumber \\
        &= \{x \in H(U) : F_1(x) \in \widetilde{D}\} \nonumber \\
        &= \{H(w) : w \in U \text{ and } \widetilde{F}(w) \in \widetilde{D}\}. \label{eq:C-repre}
    \end{align}
    Let
    \begin{equation*}
        \widetilde{C} = \{w \in U : \widetilde{F}(w) \in \widetilde{D}\}.
    \end{equation*}
    Note that $\widetilde{D}$ is a closed convex set with nonempty interior. Next we verify the constraint qualification condition. Suppose
    \begin{equation}\label{eq:CQ2}
        y \in N_{\Tilde{D}} (\widetilde{F}(0)) ~\text{ and }~ \nabla \widetilde{F} (0)^{\top} y = 0.
    \end{equation}
    Then we have
    \begin{equation*}
        \nabla F_1(\Bar{x})^{\top} y \in \mathrm{Ker}(\nabla H(0)^{\top}) = N_{F^{-1}(K)} (\Bar{x}) = \mathrm{Range}(\nabla F_2(\Bar{x})^{\top}).
    \end{equation*}
    Hence there exists $z \in \R^{m-k}$ such that
    \begin{equation*}
        \nabla F(\Bar{x})^{\top} (y, z) = \nabla F_1(\Bar{x})^{\top} y + \nabla F_2(\Bar{x})^{\top} z = 0.
    \end{equation*}
    Moreover, it is easy to deduce from (\ref{eq:CQ2}) that $(y, z) \in N_D (F(\Bar{x}))$. So by (\ref{eq:CQ}) we have $y = 0$. Thus, $\widetilde{C}$ is amenable at $0$. Applying Proposition \ref{prop:amen-int}, we know that $\widetilde{C}$ is smoothly approximately convex at $0$. Consequently, by (\ref{eq:C-repre}) and Lemma \ref{lem:sac-img}, $C$ is smoothly approximately convex at $\Bar{x}$.

\finpf
\begin{rem}
    Obviously, amenability is persistent nearby. Hence, combining Remark \ref{rmk:sac-nonlocal} and Theorem \ref{thm:euclidean}, we can see that smooth approximate convexity does not necessarily imply amenability. Specifically, the set $\mathrm{epi}\, f$, where $f$ is defined as in Example \ref{ex:super-sac}, is smoothly approximately convex at the point $(0, 0)$, but not amenable there. (Otherwise, we would have amenability and therefore smooth approximate convexity around $(0, 0)$ as well.)
\end{rem}

\subsection{Prox-regularity vs. smooth approximate convexity}
\label{sec:ecli-prox}

Next we show that, likewise to amenability, prox-regularity is a strictly stronger condition than smooth approximate convexity. Before proving this, we note that there is no implication relationship between amenability and prox-regularity. We present two examples below to illustrate this. The first example defines a set that is everywhere amenable but not prox-regular at $(0, 0)$. The second example describes a set that is prox-regular at $(0, 0)$, but not amenable there.

\begin{exa}[Amenability does not imply prox-regularity]
\label{ex:amen-noprox}
    Consider the set 
    \[C = \{(x, |x|^{3/2}) : x \in \R\}.\] 
    Clearly, $C = \{ (x, y) \in \R^2 : F(x, y) = 0\}$, where $F(x, y) = |x|^{3/2} - y$ is a $\C^{(1)}$-smooth mapping with $\nabla F$ everywhere surjective. Therefore, $C$ is everywhere amenable. On the other hand, given any $r > 0$, for all $x \in \R$ close to $0$ we have
    \begin{equation*}
        \langle (0, 1), (x, |x|^{3/2}) - (0, 0) \rangle > r \|(x, |x|^{3/2}) - (0, 0)\|^2.
    \end{equation*}
    So $C$ is not prox-regular at $(0, 0)$ for $(0, 1) \in N_C(0, 0)$.
\end{exa}

\begin{exa}[Prox-regularity does not imply amenability]
    Consider the set 
    \[C = \{(x, y) \in \R^2 : -x^2 \leq y \leq x^2\}.\] 
    It is easy to check by definition that $C$ is prox-regular at $(0, 0)$. For the sake of contradiction, assume $C$ is amenable there with representation $(V, F, D)$. Without loss of generality assume $F(0, 0) = \mathbf{0}$, where $\mathbf{0}$ is the zero vector in $\R^m$.
    
    By \cite[Theorem 6.14]{VA}, 
    \begin{equation*}
        \nabla F(0, 0)^{\top} N_D (\mathbf{0})  = N_C (0, 0) = \{0\} \times \R,
    \end{equation*}
    which, along with the constraint qualification, implies that $N_D (\mathbf{0})$ is a $1$-dimensional subspace of $\R^m$. We deduce that $\mathbf{0} \in \mathrm{ri}(D)$. Otherwise, by a standard separation argument we can find a nonzero vector $v \in \mathrm{aff}(D)$ such that $v \in N_D (\mathbf{0}) = \mathrm{aff}(D)^{\perp}$, which is a contradiction.
    
    Thus we have $C = F^{-1} (\mathrm{aff}(D))$ around $(0, 0)$, which along with the constraint qualification implies that $C$ is a $\C^{(1)}$ manifold around $(0, 0)$. However, this is obviously not the case. Hence the contradiction.
\end{exa}

\begin{rem}
    It is important that the mapping $F$ in Example \ref{ex:amen-noprox} is only $\C^{(1)}$-smooth. As explained in the Introduction, when the mapping $F$ in the representation of an amenable set is $\C^{(2)}$-smooth, the set is called ``strongly amenable''. Strongly amenable sets are always prox-regular. We can easily check this by using the second-order Taylor expansion of $F$ and the expression $N_C (\Bar{x}) = \nabla F(\bar{x})^{\top} N_D (F(\Bar{x}))$ (see \cite{VA}).
\end{rem}

Here comes the second main result. In the absence of direct implication relations, the following result, along with Theorem \ref{thm:euclidean}, helps us draw some parallel between prox-regularity and amenability.

\begin{thm}[Prox-regularity implies smooth approximate convexity]
\label{thm:prox-euclidean}
    If a set $C \subseteq \R^n$ is prox-regular at a point $\bar{x} \in C$, then $C$ is smoothly approximately convex at $\bar{x}$. Moreover, for any $\epsilon > 0$, around $\bar x$, every averaging map in $(C, d_C)$ is an $\epsilon$-path.
\end{thm}

\pf
    By \cite[Proposition 3.15 and Remark 3.2]{lewis-lopez-nicolae-22}, there exist $\sigma, r_0 > 0$ such that for all $r \in (0, r_0)$, the set $C \cap \overline{B}_r (\bar{x})$ is closed with
    \begin{equation}\label{eq:fin-ex}
        d_{C \cap \overline{B}_r (\bar{x})} (y, z) \leq \|y - z\| (1 + \sigma \|y - z\|^2) \qquad \mbox{for all}~  y, z \in C \cap \overline{B}_r (\bar{x}).
    \end{equation}
    This condition is discussed in more detail in Appendix \ref{sec:equi}.
    
    Fix any $\epsilon \in (0, 1)$. Take $0 < r < \min\left\{ r_0, \frac{\epsilon}{260 \sqrt{\sigma}} \right\}$. Denote $D = C \cap \overline{B}_r (\bar{x})$. Consider $x, x' \in D$ and let $\gamma: [0, 1] \rightarrow D$ be an averaging map in $(D, d_D)$ from $x$ to $x'$, namely a curve satisfying
    \begin{equation*}
        d_D \left( \gamma(t_1), \gamma(t_2) \right) = |t_1 - t_2| \cdot d_D (x, x') \qquad \mbox{for all}~  t_1, t_2 \in [0, 1],
    \end{equation*}
    which exists since $(D, d_D)$ is a geodesic space. We will show that $\gamma$ is an $\epsilon$-path from $x$ to $x'$.
    
    Fix any $t_1, t_2 \in [0, 1]$. Denote $m = \frac{1}{2} \gamma (t_1) + \frac{1}{2} \gamma (t_2)$ and $m' = \gamma \left( \frac{t_1 + t_2}{2} \right)$. Then
    \begin{align*}
        \|m' - m\|^2
        &= \frac{1}{2} \|m' - \gamma(t_1)\|^2 + \frac{1}{2} \|m' - \gamma(t_2)\|^2 - \frac{1}{4} \|\gamma(t_1) - \gamma(t_2)\|^2 \\
        &\leq \frac{1}{2} d_D \left(m', \gamma (t_1)\right)^2 + \frac{1}{2} d_D \left(m', \gamma (t_2)\right)^2 - \frac{1}{4} \|\gamma(t_1) - \gamma(t_2)\|^2 \\
        &= \frac{1}{4} d_D \left(\gamma (t_1), \gamma (t_2)\right)^2 - \frac{1}{4} \|\gamma(t_1) - \gamma(t_2)\|^2 \\
        &\leq \frac{\sigma}{4} \|\gamma(t_1) - \gamma(t_2)\|^4 \left( 2 + \sigma \|\gamma(t_1) - \gamma(t_2)\|^2 \right),
    \end{align*}
    where the last inequality uses (\ref{eq:fin-ex}). Since $\sigma \|\gamma(t_1) - \gamma(t_2)\|^2 \leq \sigma (2r)^2 < 1$ and
    \begin{equation*}
        \|\gamma(t_1) - \gamma(t_2)\| \leq d_D \left( \gamma (t_1), \gamma (t_2) \right) \leq |t_1 - t_2| \cdot 2 \|x - x'\|,
    \end{equation*}
    we have
    \begin{equation}\label{eq:mm}
        \|m' - m\| \leq \sqrt{\sigma} \|\gamma(t_1) - \gamma(t_2)\|^2 \leq 4 \sqrt{\sigma} \|x - x'\|^2 |t_1 - t_2|^2 \leq K |t_1 - t_2|^2,
    \end{equation}
    where $K := 8 \sqrt{\sigma} r \|x - x'\|$.
    
    For all $t, h \in [0, 1]$ such that $t + 2h \in [0, 1]$, plugging $t_1 = t, t_2 = t + 2h$ into (\ref{eq:mm}) yields
    \begin{equation*}
        \| \gamma(t) - 2 \gamma (t + h) + \gamma (t + 2h) \| \leq 8 K h^2.
    \end{equation*}
    Hence, by \cite[Lemma 1]{Philip58}, $\gamma$ is a smooth curve segment in $C$.
    
    It remains to show that $\gamma$ satisfies (\ref{eq:def-sac}). To see this, define the function $g: \overline{B}_r (\bar{x}) \rightarrow \R$ by $g(y) = d\left( y, [x, x'] \right)$ and let $f = g \circ \gamma : [0, 1] \rightarrow \R$. By the $1$-Lipschitzness and convexity of $g$, we have
    \begin{align*}
        f \left( \frac{t_1 + t_2}{2} \right) &= g(m') \leq g(m) + \|m' - m\| \leq \frac{1}{2} f(t_1) + \frac{1}{2} f(t_2) + \|m' - m\| \\
        &\leq \frac{1}{2} f(t_1) + \frac{1}{2} f(t_2) + K |t_1 - t_2|^2\qquad \mbox{for all}~ t_1, t_2 \in [0, 1].
    \end{align*}
    Using a standard argument, we deduce that
    \begin{equation*}
        f \left( (1-\alpha) t_1 + \alpha t_2 \right) \leq (1-\alpha) f(t_1) + \alpha f(t_2) + 4 K \alpha (1-\alpha) |t_1 - t_2|^2 \qquad \mbox{for all}~\alpha, t_1, t_2 \in [0, 1].
    \end{equation*}
    
    Now fix any $t \in [0, 1]$. Plugging $\alpha = t, t_1 = 0, t_2 = 1$ into the above inequality, we obtain
    \begin{equation*}
        f(t) \leq 4 K t (1-t) \leq 32 \sqrt{\sigma} r t \|x' - x\|.
    \end{equation*}
    Let $p(t) = P_{[x, x']} (\gamma (t))$. Then
    \begin{equation*}
        \|p(t) - ((1-t) x + t x')\| = \left| \|p(t) - x\| - t \|x' - x\| \right| \leq \left| \|\gamma(t) - x\| - t \|x' - x\| \right| + f(t).
    \end{equation*}
    Note that
    \begin{equation*}
        \|\gamma(t) - x\| \leq d_D (x, \gamma(t)) = t d_D (x, x') \leq (1 + 4 \sigma r^2) t\|x' - x\|
    \end{equation*}
    and
    \begin{equation*}
        \|\gamma(t) - x\| \geq \frac{d_D (x, \gamma(t))}{1 + 4 \sigma r^2} = \frac{t d_D (x, x')}{1 + 4 \sigma r^2} \geq \frac{t\|x' - x\|}{1 + 4 \sigma r^2},
    \end{equation*}
    which means
    \begin{equation*}
        \left| \|\gamma(t) - x\| - t \|x' - x\| \right| \leq 4 \sigma r^2 t\|x' - x\| \leq \sqrt{\sigma} r t\|x' - x\|.
    \end{equation*}
    Hence,
    \begin{align}
        \|\gamma(t) - ((1-t) x + t x')\| 
        &\leq \left| \|\gamma(t) - x\| - t \|x' - x\| \right| + 2 f(t) \nonumber \\
        &\leq 65 \sqrt{\sigma} r t \|x' - x\| \leq \frac{\epsilon}{4} t \|x' - x\|. \nonumber
    \end{align}
    Therefore,
    \begin{equation} \label{eq:prox-sac-1}
        \|\gamma'(0) - (x'-x)\| \le \frac{\epsilon}{4} \|x'-x\|.
    \end{equation}

    Finally, for any $t_0 \in (0, 1/2]$, let $x_0 = \gamma(t_0)$ and define the curve
    \begin{equation*}
        \Tilde{\gamma}(t) = \gamma \big(t_0 + t(1-t_0)\big)\qquad \mbox{for all}~ t \in [0, 1].
    \end{equation*}
    Then for all $t_1, t_2 \in [0, 1]$,
    \begin{equation*}
        d_D \left( \Tilde{\gamma}(t_1), \Tilde{\gamma}(t_2) \right) = (1-t_0)|t_1 - t_2| d_D (x, x') = |t_1 - t_2| d_D (x_0, x').
    \end{equation*}
    Applying conclusion \eqref{eq:prox-sac-1} to the curve $\Tilde{\gamma}$, we obtain
     \begin{align*}
        \|\Tilde{\gamma}'(0) - (x'-x_0)\|
        &\le \frac{\epsilon}{4} \|x'-x_0\| \le \frac{\epsilon}{4} d_D (x', x_0) \\
        &= \frac{\epsilon}{4} (1-t_0) d_D (x', x) \le \frac{\epsilon}{2} (1-t_0) \|x'-x\|.
    \end{align*}
    Hence,
    \begin{align*}
        \|\gamma'(t_0) - (x'-x)\|
        &= \frac{1}{1-t_0} \|\Tilde{\gamma}'(0) - (1-t_0)(x'-x)\| \\
        &\le \frac{1}{1-t_0} \|\Tilde{\gamma}'(0) - (x'-x_0)\| + \frac{1}{1-t_0} \|(1-t_0)x + t_0x' - x_0\| \\
        &\le \frac{\epsilon}{2} \|x'-x\| + \frac{t_0}{4(1-t_0)} \epsilon \|x'-x\| \le \epsilon \|x'-x\|.
    \end{align*}
    By symmetry, we know that the above inequality also holds for $t_0 \in [1/2, 1]$.
    
    The second statement follows by noting that all points $y, y'$ sufficiently close to $\bar x$ satisfy $d_D (y, y') = d_C (y, y')$. Thus every averaging map in $(C, d_C)$ sufficiently close to $\bar x$ is also an averaging map in $(D, d_D)$.
\finpf

\begin{rem}
    It is clear from Example \ref{ex:amen-noprox} that the converse to the above implication does not hold.
\end{rem}

\section{Smooth approximate convexity on Riemannian manifolds}
\label{sec:manifold}

In light of Lemma \ref{lem:sac-img}, we naturally realize that the structure of smooth approximate convexity is not limited to Euclidean spaces. In this section, we generalize the idea and the two main results to Riemannian manifolds. We first consider embedded submanifolds in Euclidean spaces and then use the Nash embedding theorem \cite{Nash56} to deal with general Riemannian manifolds.

\subsection{Embedded submanifolds}\label{sec:embedded}

The generalization to embedded submanifolds is straightforward. In Euclidean spaces, smooth approximate convexity is characterized by comparing the derivatives of the paths with those of line segments. In embedded manifolds, we simply replace the line segments with averaging maps (see \eqref{eq:def-averaging-map}).

We consider a $\C^{(3)}$-smooth embedded submanifold $\M$ in $\R^n$ and a point $\bar x \in \M$. In this case, every pair of points in $\M$ near $\bar x$ are connected by a unique geodesic, therefore also by a unique averaging map. Moreover, since $\M$ is prox-regular at $\bar x$ (see \cite{VA}), we can apply Theorem \ref{thm:prox-euclidean} to see that all averaging maps of $\M$ near $\bar x$ are smooth curve segments. Note that $\C^{(3)}$-smoothness of the manifold is assumed for simplicity. Most of the results below hold under weaker smoothness conditions.

\begin{defn}\label{def:sac-rm}
    A set $C \subseteq \mathcal{M}$ is {\em smoothly approximately convex} at a point $\Bar{x} \in C$ in $\M$ if given any $\epsilon > 0$, for all distinct $x, x' \in C$ around $\Bar{x}$, there is a smooth curve segment $\gamma: [0, 1] \rightarrow C$ with $\gamma(0) = x$ and $\gamma(1) = x'$ that satisfies
    \begin{equation}\label{eq:def-sac-rm}
        \|\gamma'(t)- g'(t)\| \leq \epsilon r \qquad \text{for all}~ t \in [0, 1],
    \end{equation}
    where $r = d_{\M}(x, x')$ and $g: [0, 1] \rightarrow \mathcal{M}$ is the averaging map from $x$ to $x'$. We call such a map $\gamma$ an {\em $\epsilon$-path} from $x$ to $x'$ in $\M$.
\end{defn}

The ideas of amenability and prox-regularity also naturally generalize.

\begin{defn}\label{def:amenable-rm}
    A set $C \subseteq \mathcal{M}$ is {\em amenable} at $\Bar{x} \in C$ in $\M$ if there is an open neighborhood $V$ of $\Bar{x}$ in $\mathcal{M}$ along with a $\C^{(1)}$ mapping $F$ from $V$ into a space $\R^m$ and a closed, convex set $D\subseteq \R^m$ such that
    \begin{equation*}
        C \cap V = \left\{x \in V : F(x) \in D \right\}
    \end{equation*}
    and the following constraint qualification is satisfied:
    \begin{equation}\label{eq:CQ-rm}
      \text{the only vector}~ y \in N_D (F(\Bar{x})) ~\text{with}~ y \perp \mathrm{Range}(d F (\bar{x})) ~\text{is}~ y=0.
    \end{equation}
    In this case, we call $(V, F, D)$ a {\em representation} of the set $C$ in $\mathcal{M}$ around $\Bar{x}$.
\end{defn}

\begin{defn}\label{def:prox-rm}
    A set $C \subseteq \mathcal{M}$ is {\em prox-regular} at $\Bar{x} \in C$ in $\M$ if there exists $r > 0$ such that $P_C (y)$ is a singleton for all $y \in B_r (\Bar{x})$.
\end{defn}

The notion of prox-regularity in Riemannian manifolds was introduced in \cite{hosseini-pour} with a different definition involving the exponential mapping. However, prox-regular sets in the sense of \cite{hosseini-pour} also satisfy the projection property in Definition \ref{def:prox-rm}
(see \cite[Theorem 3.12]{hosseini-pour}). Note also that Definition \ref{def:prox-rm} implies that $C$ is locally closed at $\bar x$.

Before proving the implications, we show that all of the generalized definitions above are compatible with the canonical definitions given in Section \ref{sec:eclidean}. Namely, whether a set satisfies the desired conditions does not depend on the ambient space that we view it in. To be specific, we have the following equivalence results.

\begin{lem}\label{lem:sac-rm}
    Given a set $C \subseteq \mathcal{M}$ and a point $\bar{x} \in C$, $C$ is smoothly approximately convex at $\bar{x}$ in $\R^n$ if and only if it is smoothly approximately convex at $\bar{x}$ in $\mathcal{M}$.
\end{lem}

\pf 
    We have seen in Theorem \ref{thm:prox-euclidean} that for all $x, x' \in \mathcal{M}$ around $\bar{x}$,
    \begin{equation}\label{eq:dm}
        d_\mathcal{M} (x, x') = \|x-x'\| + O(\|x-x'\|^3).
    \end{equation}
Therefore to prove the desired result, it is enough to show that given any $\epsilon > 0$, for all distinct $x, x' \in \mathcal{M}$ around $\Bar{x}$, the averaging map $g: [0, 1] \rightarrow \mathcal{M}$ from $x$ to $x'$ satisfies 
    \begin{equation*}
        \|g'(t)- (x'-x)\| \leq \epsilon \|x'-x\| \qquad \text{for all}~ t \in [0, 1],
    \end{equation*}
    namely that $g$ is an $\epsilon$-path in $\R^n$. This holds by Theorem \ref{thm:prox-euclidean}.
\finpf

Similar correspondences hold for amenability and prox-regularity. 

\begin{lem}\label{lem:amenable}
    Given a set $C \subseteq \mathcal{M}$ and a point $\bar{x} \in C$, $C$ is amenable at $\bar{x}$ in $\R^n$ if and only if it is amenable at $\bar{x}$ in $\mathcal{M}$.
\end{lem}

\begin{lem}\label{lem:prox-regularity}
    Given a set $C \subseteq \mathcal{M}$ and a point $\bar{x} \in C$, $C$ is prox-regular at $\bar{x}$ in $\R^n$ if and only if it is prox-regular at $\bar{x}$ in $\mathcal{M}$.
\end{lem}

The proof for the former involves illuminating results on the structure of amenable sets. In connection to the latter, Bangert \cite{Ban82} studied sets that are prox-regular at all their points and showed, based on \cite{Kleinjohann81}, that this notion does not depend on the considered Riemannian metric. Instead of building a proof relying on this result, we justify Lemma \ref{lem:prox-regularity} as a consequence of the fact that prox-regularity at some point implies finite extrinsic curvature there. To avoid distracting readers from the main focus of the paper, we discuss these proofs in the appendices.

Armed with the above lemmas, we readily obtain the following generalized results.

\begin{thm}\label{thm:embedded}
    If a set $C$ is amenable at a point $\bar{x} \in C$ in $\mathcal{M}$, then it is smoothly approximately convex at $\bar{x}$ in $\mathcal{M}$.
\end{thm}

\pf
    Follows directly from Lemma \ref{lem:sac-rm}, Lemma \ref{lem:amenable}, and Theorem \ref{thm:euclidean}.
\finpf

\begin{thm}\label{thm:prox-embedded}
    If a set $C$ is prox-regular at a point $\bar{x} \in C$ in $\mathcal{M}$, then it is smoothly approximately convex at $\Bar{x}$ in $\mathcal{M}$.
\end{thm}

\pf
    Follows directly from Lemma \ref{lem:sac-rm}, Lemma \ref{lem:prox-regularity}, and Theorem \ref{thm:prox-euclidean}.
\finpf

\subsection{General manifolds}
\label{sec:rm}

Next we move on to manifolds that are not necessarily subsets of Euclidean spaces. Readers unfamiliar with general Riemannian manifolds can refer to \cite{Lee, lee-riemannian}. Assume $\mathcal{M}$ is a smooth Riemannian manifold, namely it is a $\C^{(\infty)}$-smooth manifold endowed with a $\C^{(\infty)}$-smooth Riemannian metric. Thanks to the Nash embedding theorem \cite{Nash56}, we can study geometry in $\M$ through an embedded submanifold, and hence establish results similar to Theorems \ref{thm:embedded} and \ref{thm:prox-embedded} in this general setting.

Let $d_\mathcal{M} (\cdot, \cdot)$ be the Riemannian distance function on $\mathcal{M}$. Given a set $C \subseteq \mathcal{M}$, we call a curve $\gamma: [a, b] \rightarrow C$ a \textit{smooth curve segment in $C$} if it has an extension to a $\C^{(1)}$-smooth map from a neighborhood of $[a, b]$ to $\mathcal{M}$. Note that in the special case where $\mathcal{M}$ is an embedded submanifold, the definitions of $d_\mathcal{M}(\cdot, \cdot)$ and of smooth curve segments here coincide with the ones in the previous sections.

We define amenability and prox-regularity in exactly the same way as in the embedded case, namely by Definitions \ref{def:amenable-rm} and \ref{def:prox-rm}. However, the concept of smooth approximate convexity does not have a direct generalization to the new setting. The reason is that, when $\M$ is not a subset of a Euclidean space, tangent vectors to $\M$ at different base points belong to different tangent spaces and are therefore impossible to compare directly. As a workaround, we consider the following formally weaker definition that retains the geometric feature of smooth approximate convexity.

\begin{defn}\label{def:sac-rm-weak}
    A set $C \subseteq \mathcal{M}$ is {\em smoothly weakly convex} at a point $\Bar{x} \in C$ in $\M$ if given any $\epsilon > 0$, for all distinct $x, x' \in C$ around $\Bar{x}$, there is a smooth curve segment $\gamma: [0, 1] \rightarrow C$ that satisfies
    \begin{equation}\label{eq:sac-dist-rm}
        d_\mathcal{M} (\gamma(t), g(t)) \leq \epsilon r t (1 - t)\qquad \mbox{for all}~ t \in [0, 1],
    \end{equation}
    and
    \begin{equation}\label{eq:sac-len-rm}
        L(\gamma) \le (1+\epsilon) r,
    \end{equation}
    where $r = d_{\M}(x, x')$ and $g: [0, 1] \rightarrow \mathcal{M}$ is the averaging map from $x$ to $x'$ (see \eqref{eq:def-averaging-map}). We call such a map $\gamma$ a {\em weak $\epsilon$-path} from $x$ to $x'$ in $\M$.
\end{defn}

The first lemma shows that in embedded submanifolds, smooth approximate convexity is a stronger condition than smooth weak convexity.

\begin{lem}\label{lem:sac-int}
    Suppose $\M$ is a $\C^{(3)}$-smooth embedded submanifold in $\R^n$. Consider any set $C \subseteq \M$, point $\bar x \in C$, and constant $\epsilon > 0$. Then around $\bar x$, every $\epsilon$-path in $\M$ is a weak $3\epsilon$-path in $\M$.
\end{lem}

\pf
    Consider points $x, x' \in \M$ around $\bar x$. Let $\gamma: [0, 1] \rightarrow C$ be an $\epsilon$-path from $x$ to $x'$, and let $g: [0, 1] \rightarrow \M$ be the averaging map. The inequality (\ref{eq:sac-len-rm}) follows easily from (\ref{eq:def-sac-rm}) and the fact that $\|g'(t)\|=r$ for all $t \in (0, 1)$, where $r = d_{\M} (x, x')$. By H\"older's inequality, for all $t \in [0, 1]$,
    \begin{align*}
        \|\gamma(t) - g(t)\| 
        &= \| \int_0^t \left(\gamma' (\tau) - g'(\tau)\right) \mathrm{d}\tau \| \\
        &\leq \sqrt{t \int_0^t \|\gamma' (\tau) - g'(\tau)\|^2 \mathrm{d}\tau} \\
        &\leq \epsilon r t.
    \end{align*}
    By the symmetry of $\gamma$ with respect to $x$ and $x'$, we also have
    \begin{equation*}
        \|\gamma(t) - g(t)\| \leq \epsilon r (1-t) \qquad \mbox{for all}~ t \in [0,1].
    \end{equation*}
    Hence,
    \begin{equation*}
        \|\gamma (t) - g(t) \| \leq 2\epsilon r t(1-t) \qquad \mbox{for all}~ t \in [0,1].
    \end{equation*}
    We conclude by (\ref{eq:dm}) that $\gamma$ is a weak $3\epsilon$-path.
\finpf

Although Definition \ref{def:sac-rm-weak} is formally weaker than smooth approximate convexity, the two inequalities in it capture the key idea of approximating geodesics with smooth feasible paths. Inequality (\ref{eq:sac-dist-rm}), with roots in the earlier ideas of approximate convexity and the UAG property (see (\ref{eq:def-intrinsic}), (\ref{eq:def-fac}), and (\ref{eq:def-uag})), measures the distance between the path and the reparametrized geodesic. Inequality (\ref{eq:sac-len-rm}) ensures that the path does not have erratic curvature, and in particular guarantees the set to be normally embedded. The weakening allows us to establish counterparts of Theorems \ref{thm:embedded} and \ref{thm:prox-embedded} in general manifolds.

As a first step, we note that smooth weak convexity possesses the same independence of the ambient space as the three conditions discussed in Section \ref{sec:embedded}.

\begin{lem}\label{lem:sac-rm-weak}
     Suppose $\M$ is a $\C^{(3)}$-smooth embedded submanifold in $\R^n$. Consider any set $C \subseteq \M$ and point $\bar x \in C$. Then $C$ is smoothly weakly convex at $\bar{x}$ in $\R^n$ if and only if it is smoothly weakly convex at $\bar{x}$ in $\mathcal{M}$.
\end{lem}

\pf In light of (\ref{eq:dm}), it suffices to prove that given any $\epsilon > 0$, for all distinct $x, x' \in \mathcal{M}$ around $\Bar{x}$, the averaging map $g: [0, 1] \rightarrow \mathcal{M}$ from $x$ to $x'$ satisfies
    \begin{equation*}
       \big \|g(t)- ((1-t)x+tx')\big\| \leq \epsilon \|x'-x\| t(1-t) \qquad \text{for all}~ t \in (0, 1),
    \end{equation*}
    namely that $g$ is a weak $\epsilon$-path in $\R^n$. This holds by Theorem \ref{thm:prox-euclidean} and Lemma \ref{lem:sac-int}.
\finpf

The next two results generalize Theorems \ref{thm:embedded} and \ref{thm:prox-embedded}.

\begin{thm}\label{thm:rm}
    If a set $C$ is amenable at $\bar{x} \in C$ in $\mathcal{M}$, then it is smoothly weakly convex there.
\end{thm}

\pf
    Using the Nash embedding theorem \cite{Nash56}, we know there is a $\C^{(3)}$ isometric embedding $H$ from $\mathcal{M}$ to a space $\R^n$. Let $\bar{w} = H(\bar{x})$, $\widetilde{\mathcal{M}} = H(\mathcal{M})$ and $\widetilde{C} = H(C)$. Then $\widetilde{\mathcal{M}}$ is a $\C^{(3)}$-smooth embedded submanifold of $\R^n$. Let $H^{-1}$ denote the inverse of the mapping $H: \mathcal{M} \rightarrow \widetilde{\mathcal{M}}$. Suppose $C \cap V = \{x \in V : F(x) \in D\}$, where $V$ is a neighborhood of $\Bar{x}$ in $\mathcal{M}$, $F: V \rightarrow \R^m$ is a $\C^{(1)}$-smooth map and $D \subseteq \R^m$ is a closed convex set satisfying (\ref{eq:CQ-rm}). Define $\widetilde{V} = H(V)$ and $\widetilde{F} = F \circ H^{-1}: \widetilde{V} \rightarrow \R^m$. Note that $d \widetilde{F}(\bar{w}) = d F(\bar{x}) \circ d H^{-1} (\bar{w})$. So $\widetilde{V}$ is a neighborhood of $\bar{w}$ in $\widetilde{M}$ satisfying
    \begin{equation*}
        \widetilde{C} \cap \widetilde{V} = \{ H(x) : x \in V, F(x) \in D \} = \{ w \in \widetilde{V} : \widetilde{F} (w) \in D \}
    \end{equation*}
    and the property that the only vector $y \in N_D (\widetilde{F}(\bar{w}))$ with $y \perp \mathrm{Range} (d \widetilde{F}(\bar{w}))$ is $y=0$. This means that $\widetilde{C}$ is amenable at $\bar{w}$ in $\widetilde{\mathcal{M}}$. By Theorem \ref{thm:embedded}, we know that $\widetilde{C}$ is smoothly approximately convex and therefore smoothly weakly convex at $\bar{w}$ in $\widetilde{\mathcal{M}}$. Fix any $\epsilon > 0$. For all distinct $x, x' \in C$ around $\bar{x}$, the points $w := H(x), w' := H(x')$ are in $\widetilde{C}$ and close to $\bar{w}$. So there is a smooth curve segment $\Tilde{\gamma}: [0, 1] \rightarrow \widetilde{C}$ that satisfies
    \begin{equation*}
        d_{\widetilde{\mathcal{M}}} (\Tilde{\gamma}(t), \Tilde{g}(t)) \leq \epsilon r t (1 - t)\qquad \text{for all}~ t \in [0, 1],
    \end{equation*}
    and
    \begin{equation*}
        L(\Tilde{\gamma}) \le (1+\epsilon) r,
    \end{equation*}
    where $r = d_{\widetilde{\mathcal{M}}} (w, w')$ and $\Tilde{g}: [0, 1] \rightarrow \widetilde{\mathcal{M}}$ is the averaging map from $w$ to $w'$. Since $H$ is an isometry, it keeps the lengths of curves and Riemannian distances around $\bar{x}$. Hence, letting $\gamma = H^{-1} \circ \Tilde{\gamma}$ and $g = H^{-1} \circ \Tilde{g}$, we know $\gamma([0, 1]) \subseteq C$, $r = d_{\M} (x, x')$, and $g$ is the averaging map from $x$ to $x'$ satisfying
    \begin{equation*}
        d_\mathcal{M} (\gamma(t), g(t)) = d_{\widetilde{\mathcal{M}}} (\Tilde{\gamma}(t), \Tilde{g}(rt)) \leq \epsilon r t (1 - t)\qquad \text{for all}~ t \in [0, 1]
    \end{equation*}
    and
    \begin{equation*}
        L(\Tilde{\gamma}) = L(\gamma) \le (1+\epsilon) r.
    \end{equation*}
    We conclude that $C$ is smoothly weakly convex at $\bar{x}$ in $\M$.
\finpf

It is clear from the definition that prox-regularity is also preserved under isometric embeddings. Hence we have the following result.

\begin{thm}\label{thm:prox-rm}
    If a set $C$ is prox-regular at $\bar{x} \in C$ in $\mathcal{M}$, then it is smoothly weakly convex there.
\end{thm}

\begin{rem}
    Although it is common practice to assume $\C^{(\infty)}$-smoothness, the Nash embedding theorem only requires $\C^{(3)}$-smoothness of the metric for the mapping $H$ to exist. Hence, Theorems \ref{thm:rm} and \ref{thm:prox-rm} also hold for any $\C^{(\infty)}$-smooth manifold endowed with a $\C^{(3)}$-smooth Riemannian metric.
\end{rem}

\section{First-order optimality condition on Riemannian manifolds}
\label{sec:optimality}

In this section, we consider the optimization problem
\begin{equation}\label{eq:min}
    \min_{x \in C}\, f(x),
\end{equation}
where $C$ is a subset of a Riemannian manifold $\mathcal{M}$, and $f: \mathcal{M} \rightarrow \R$ is a $\C^{(1)}$-smooth function. To discuss the first-order optimality condition, we need to introduce the following property for a tangent vector $v$ at $x \in C$, which is equivalent to the fact that $v$ belongs to the Bouligand tangent cone to $C$ at $x$ as given in \cite[Definition 3.8]{ledyaev-zhu}. 

\begin{defn}
	Let $x \in C$. A vector $v \in T_\mathcal{M} (x)$ is {\em tangent to $C$ at $x$}, written $v \in T_C (x)$, if $v = 0$ or there exists a sequence $\{x_k\} \subseteq C$ converging to $x$ such that $g_k'(0) \rightarrow v/\|v\|$, where $g_k$ is the geodesic from $x$ to $x_k$.
\end{defn}

Consider the following standard first-order optimality condition at a point $\bar{x} \in C$:
\begin{equation}\label{eq:1st-order}
    v f \geq 0 \qquad \text{for all}~v \in T_C (\bar{x}).
\end{equation}
This condition must hold at any local minimizer of $f$ on $C$. The proof is elementary (cf. \cite[Proposition 3.9]{ledyaev-zhu}). 

An interesting property of smoothly approximately convex sets is that for such sets, the failure of (\ref{eq:1st-order}) implies the existence of a $\C^{(1)}$-smooth descent path of $f$ in $C$ starting from $\bar{x}$. 

\begin{prop}\label{prop:opt-rm}
    Let $\mathcal{M}$ be a smooth Riemannian manifold. If a set $C$ is smoothly weakly convex at a point $\bar{x} \in C$ in $\mathcal{M}$, then the failure of (\ref{eq:1st-order}) implies that there is a smooth curve segment in $C$ starting from $\bar{x}$ on which $f$ is decreasing at a linear rate.
\end{prop}

\pf
    By the Nash embedding theorem \cite{Nash56}, we know there is a $\C^{(3)}$ isometric embedding $H$ from $\mathcal{M}$ to a space $\R^n$. Let $\bar{w} = H(\bar{x})$, $\widetilde{\mathcal{M}} = H(\mathcal{M})$ and $\widetilde{C} = H(C)$. Then $\widetilde{\mathcal{M}}$ is a $\C^{(3)}$-smooth embedded submanifold of $\R^n$. Let $H^{-1}$ denote the inverse of the mapping $H: \mathcal{M} \rightarrow \widetilde{\mathcal{M}}$ and define $\Tilde{f} = f \circ H^{-1}$. Using the same argument as in Theorem \ref{thm:rm}, we can see that $\widetilde{C} \subseteq \widetilde{\mathcal{M}}$ is smoothly weakly convex at $\bar{w}$. Suppose there exists a vector $v \in T_C (\bar{x})$ such that $v f < 0$. It is easy to check that
    \begin{equation*}
        \Tilde{v} := d H(\bar{x}) v \,\in\, T_{\widetilde{C}} (\bar{w}) \quad \textup{and} \quad \Tilde{v} \Tilde{f} = v f < 0.
    \end{equation*}    
    By definition there is a sequence $\{w_k\} \subseteq \widetilde{C}$ converging to $\bar{w}$ such that $g_k'(0) \rightarrow \Tilde{v}$, where $g_k: [0, r_k] \to \widetilde{\mathcal{M}}$ is the geodesic from $\bar{w}$ to $w_k$. Let $\bar{f}$ be a $\C^{(1)}$-smooth function defined on a neighborhood of $\widetilde{\mathcal{M}}$ that satisfies $\bar{f}|_{\widetilde{\mathcal{M}}} = \Tilde{f}$. Take
    \begin{equation*}
        \epsilon := - \frac{\Tilde{v} \Tilde{f}}{2 \|\nabla \bar{f} (\bar{w})\|} = - \frac{1} {2 \|\nabla \bar{f} (\bar{w})\|} \langle \nabla \bar{f} (\bar{w}), \Tilde{v} \rangle > 0.
    \end{equation*}
    By the smooth weak convexity of $\widetilde{C}$ at $\bar{w}$ and considering that $t \mapsto g_k(r_k t)$ are averaging maps (see \eqref{eq:def-averaging-map}), we know there exist $k \in \N$ and a smooth curve segment $\gamma: [0, 1] \rightarrow \widetilde{C}$ such that $\|g_k'(0) - \Tilde{v}\| < \epsilon$ and
    \begin{equation*}
        d_{\widetilde{\mathcal{M}}} (\gamma (t), g_k (r_k t)) \leq \epsilon r_k t (1 - t) \qquad \text{for all}~ t \in [0, 1].
    \end{equation*}
    It follows that
    \begin{align*}
        \| \gamma'(0) - r_k g_k'(0) \|
        &= \lim_{t \searrow 0} \frac{1}{t} \|\gamma (t) - g_k (r_k t)\| \\
        &\leq \liminf_{t \searrow 0} \frac{1}{t} d_{\widetilde{\mathcal{M}}}(\gamma (t), g_k (r_k t)) \leq \epsilon r_k.
    \end{align*}
    Hence defining $\tilde{\gamma}: [0, r_k] \rightarrow \widetilde{C}$ by $\tilde{\gamma}(t) = \gamma (t/r_k)$ we have
    \begin{align*}
       (\Tilde{f} \circ \Tilde{\gamma})'(0)
        &= \langle \nabla \bar{f} (\bar{w}), \Tilde{\gamma}'(0) \rangle \\
        &\leq \langle \nabla \bar{f} (\bar{w}), \Tilde{v} \rangle + \| \nabla \bar{f} (\bar{w}) \| \left( \| \Tilde{\gamma}'(0) - g_k'(0) \| + \|g_k'(0) - \Tilde{v}\| \right) \\
        &< \langle \nabla \bar{f} (\bar{w}), \Tilde{v} \rangle + \| \nabla \bar{f} (\bar{w}) \| \cdot 2 \epsilon = 0.
    \end{align*}
    Thus $H^{-1} \circ \Tilde{\gamma}$ is a descent curve of $f$.
\finpf

\begin{rem}
From the proof above, we can see that the conclusion of Proposition \ref{prop:opt-rm} also holds when $\M$ is a $\C^{(3)}$-smooth embedded submanifold in $\R^n$ and $C \subseteq \M$ is smoothly weakly convex (or smoothly approximately convex) at $\bar x$ in $\M$ (or, equivalently, in $\R^n$).
\end{rem}

Feasible curves of linear decrease, like that guaranteed in Proposition \ref{prop:opt-rm}, open up the possibility of constructing descent steps using a backtracking line search, like that in \cite[Section 4.5]{boumal2022intromanifolds}.

\appendix
\normalsize

\section{Appendix: amenability on Riemannian manifolds} \label{sec:equi-amenable}

In this section, we take a few steps to prove Lemma \ref{lem:amenable}. As a byproduct, the intermediate results also help uncover the intrinsic structure of amenable sets.

Still consider a $\C^{(3)}$-smooth embedded submanifold $\M$ in $\R^n$. The first result is an analogue of Lemma \ref{lem:CQ} combined with Remark \ref{rmk:connected}.

\begin{lem}\label{lem:CQ-manifold}
    In Definition \ref{def:amenable-rm}, if $\mathrm{int}(D) \neq \emptyset$, then the constraint qualification (\ref{eq:CQ-rm}) is equivalent to the following condition:
    \begin{equation*}\label{eq:CQ-rm-dual}
    	\text{there exists}~ u \in T_{\mathcal{M}} (\Bar{x}) ~\text{such that}~ F(\Bar{x}) + d F(\Bar{x}) u \in \mathrm{int}(D).
    \end{equation*}
    Moreover, when these conditions hold, $F^{-1} (\mathrm{int}(D))$ is a connected set around $\Bar{x}$.
\end{lem}

\pf
    Suppose $\mathcal{M} = H(B_{\delta} (0))$ around $\Bar{x}$, where $B_{\delta} (0) \subseteq \R^k$ and $H: B_{\delta} (0) \rightarrow \R^n$ is an embedding. Let $\widetilde{C} = H^{-1} (C)$. Then there exists an open neighborhood $W$ of $0$ in $\R^k$ such that $\widetilde{C} \cap W = \{w \in W : \widetilde{F} (w) \in D\}$, where $\widetilde{F} := F \circ H$. Note that $\nabla \widetilde{F} (0) = d F (\Bar{x}) \nabla H(0)$ and $T_{\mathcal{M}} (\Bar{x}) = \mathrm{Range} (\nabla H(0))$. Hence the equivalence follows directly from applying Lemma \ref{lem:CQ} to $(W, \widetilde{F}, D)$. The local connectedness of $F^{-1} (\mathrm{int}(D))$ comes from  applying Remark \ref{rmk:connected} to $(W, \widetilde{F}, D)$.
    
\finpf

The following lemma allows us to understand general amenable sets.

\begin{lem}\label{lem:amenable-structure}
    Consider a set $C \subseteq \R^n$ and a point $\Bar{x} \in C$. Then $C$ is amenable at $\Bar{x}$ in $\R^n$ if and only if there is a $\C^{(1)}$-smooth embedded submanifold $\mathcal{N}$ in $\R^n$ such that $C$ is amenable at $\Bar{x}$ in $\mathcal{N}$ with a representation $(V_\mathcal{N}, F_\mathcal{N}, D_\mathcal{N})$ satisfying $\mathrm{int} (D_\mathcal{N}) \neq \emptyset$.
\end{lem}

\pf
    For the direct implication, let $(V, F, D)$ be a representation of $C$ in $\R^n$ around $\Bar{x}$. Define $F_1, F_2, P$ in the same way as in Theorem \ref{thm:euclidean}. Let
    \begin{equation*}
        \mathcal{N} = F_2^{-1} (0),~ V_{\mathcal{N}} = V \cap \mathcal{N},~ F_{\mathcal{N}} = F_1 |_{\mathcal{N}},~ D_{\mathcal{N}} = P(D).
    \end{equation*}
    Then $C \subseteq \mathcal{N}$ and $C \cap V_{\mathcal{N}} = \{x \in V_{\mathcal{N}} : F_{\mathcal{N}} (x) \in D_{\mathcal{N}} \}$. Moreover, whenever $y \in N_{D_{\mathcal{N}}} (F_{\mathcal{N}} (\Bar{x}))$ and $y \perp \mathrm{Range} (d F_{\mathcal{N}} (\Bar{x})) = \nabla F_1 (\Bar{x}) T_{\mathcal{N}} (\Bar{x})$, we have
    \begin{equation*}
        \nabla F_1 (\Bar{x})^{\top} y \,\in\, N_{\mathcal{N}} (\Bar{x}) = \mathrm{Range} (\nabla F_2 (\Bar{x})^{\top}),
    \end{equation*}
    which means there exists $z$ such that $(y, z) \in \mathrm{Ker} (\nabla F (\Bar{x})^{\top})$, therefore $y = 0$.

    Conversely, suppose $C$ is amenable at $\Bar{x}$ in $\mathcal{N}$ with representation $(V_{\mathcal{N}}, F_{\mathcal{N}}, D_{\mathcal{N}})$. According to \cite{Lee}, there is an open neighborhood $V$ of $\Bar{x}$ in $\R^n$ along with two $\C^{(1)}$-smooth mappings $F_1 : V \rightarrow \R^{k_1}$ and $F_2 : V \rightarrow \R^{k_2}$ such that $F_{\mathcal{N}} = F_1 |_{\mathcal{N}}$, $V \cap \mathcal{N} = F_2^{-1} (0) \subseteq V_{\mathcal{N}}$, and $\nabla F_2 (\Bar{x})$ is surjective. Let $F = (F_1, F_2)$ and $D = D_{\mathcal{N}} \times \{0\} \subseteq \R^{k_1} \times \R^{k_2}$. Then $C \cap V = \{x \in V: F(x) \in D\}$. Let $(y, z) \in N_D (F(\bar{x})) \cap \mathrm{Ker} (\nabla F (\Bar{x})^{\top})$. Observe that
    \begin{equation*}
        \nabla F_1 (\Bar{x})^{T} y = - \nabla F_2 (\Bar{x})^{T} z \,\in\, \mathrm{Range} (\nabla F_2 (\Bar{x})^{\top}) = N_{\mathcal{N}} (\Bar{x}).
    \end{equation*}
    Hence, $y \in N_{D_{\mathcal{N}}} (F_{\mathcal{N}} (\Bar{x})) \cap \mathrm{Range} (d F_{\mathcal{N}} (\Bar{x}))^{\perp} = \{0\}$ and we get $\nabla F_2 (\Bar{x})^{T} z = 0$. It follows from the surjectivity of $\nabla F_2 (\Bar{x})$ that $z = 0$. So $(V, F, D)$ is a valid representation of $C$ in $\R^n$ around $\Bar{x}$.
\finpf

\begin{rem}\label{rmk:amenable-structure} 
We can see that in Lemma \ref{lem:amenable-structure}, the condition $\mathrm{int} (D_{\mathcal{N}}) \neq \emptyset$ is not needed for the reverse implication. Moreover, the equivalence still holds if we replace $\R^n$ with the manifold $\mathcal{M}$.
\end{rem}

\begin{lem}\label{lem:manifold-proj}
    Consider two embedded submanifolds $\mathcal{N}_1$ and $\mathcal{N}_2$ of $\R^n$, where $\mathcal{N}_1$ is $\C^{(1)}$-smooth and $\mathcal{N}_2$ is $\C^{(2)}$-smooth. Suppose that, at a point $\Bar{x} \in \mathcal{N}_1 \cap \mathcal{N}_2$, we have $T_{\mathcal{N}_1} (\Bar{x}) \subseteq T_{\mathcal{N}_2} (\Bar{x})$. Then $\mathrm{P}_{\mathcal{N}_2} (\mathcal{N}_1)$ is also a $\C^{(1)}$-smooth embedded submanifold of $\R^n$ around $\Bar{x}$.
\end{lem}

\pf
    Suppose $\mathcal{N}_1 = H(B_{\delta}(0))$ around $\Bar{x}$, where $B_{\delta}(0)$ is an open ball in $\R^k$ and $H: B_{\delta}(0) \rightarrow \R^n$ is a $\C^{(1)}$-smooth embedding. Let $\widetilde{H} = \mathrm{P}_{\mathcal{N}_2} \circ H$. Then by \cite[Lemma 2.1]{Lewis08}, $\widetilde{H}$ is $\C^{(1)}$-smooth around $\Bar{x}$ with
    \begin{equation*}
        \nabla \widetilde{H} (0) = \nabla \mathrm{P}_{\mathcal{N}_2} (\Bar{x}) \cdot \nabla H(0) = \mathrm{P}_{T_{\mathcal{N}_2} (\Bar{x})} \nabla H(0) = \nabla H(0).
    \end{equation*}
    So $\widetilde{H}$ is also a $\C^{(1)}$-smooth embedding on $B_{\delta'}(0)$ for some $\delta' > 0$.
\finpf

We can now prove the correspondence between amenability in $\R^n$ and that in $\mathcal{M}$.

\vspace{3mm}
\noindent{\bf Proof of Lemma \ref{lem:amenable}\ \ }
    The reverse implication is already proved in Lemma \ref{lem:amenable-structure}. We focus on the direct implication. By Lemma \ref{lem:amenable-structure}, there is a $\C^{(1)}$-smooth embedded submanifold $\mathcal{N}$ in $\R^n$ such that $C$ is amenable at $\bar x$ in $\mathcal{N}$ with a representation $(V_{\mathcal{N}}, F_{\mathcal{N}}, D)$ satisfying $\mathrm{int} (D) \neq \emptyset$. Then applying Lemma \ref{lem:CQ-manifold}, we know there exists a vector $u$ such that
    \begin{equation}\label{eq:u-tc}
        u \in T_{\mathcal{N}} (\Bar{x}) \quad \text{and} \quad F_{\mathcal{N}} (\Bar{x}) + d F_{\mathcal{N}} (\Bar{x}) u \in \mathrm{int} (D).
    \end{equation}
    Let $\gamma: [0, \epsilon) \rightarrow \mathcal{N}$ be a $\C^{(1)}$-smooth curve with $\gamma'(0) = u$. The convexity of $D$ along with (\ref{eq:u-tc}) implies that $F_{\mathcal{N}} (\gamma(t)) \in \mathrm{int} (D)$ for all sufficiently small $t > 0$. Note that around every such $\gamma(t)$, we have $\mathcal{N} \subseteq \mathcal{M}$. Hence,
    \begin{equation*}
        T_{\mathcal{N}} (\Bar{x}) \,=\, \lim_{k \rightarrow 
        \infty} T_{\mathcal{N}} (\gamma(1/k)) ~\subseteq~ \lim_{k \rightarrow 
        \infty} T_{\mathcal{M}} (\gamma(1/k)) \,=\, T_{\mathcal{M}} (\Bar{x}).
    \end{equation*}
    Clearly, for every $x \in C$ sufficiently close to $\Bar{x}$, (\ref{eq:u-tc}) as well as the arguments after it, still hold when we replace $\Bar{x}$ with $x$. Therefore, $C = \mathrm{cl} (C_0)$ around $\Bar{x}$, where $C_0 := F_{\mathcal{N}}^{-1} (\mathrm{int}(D))$. 
    
    Let $\widetilde{\mathcal{N}} = \mathrm{P}_{\mathcal{M}} (\mathcal{N})$. By Lemma \ref{lem:manifold-proj}, $\widetilde{\mathcal{N}}$ is a $\C^{(1)}$-smooth embedded submanifold of $\R^n$ around $\Bar{x}$ with $T_{\widetilde{\mathcal{N}}} (\Bar{x}) = T_{\mathcal{N}} (\Bar{x})$. Note that $C \subseteq \mathcal{N} \cap \widetilde{\mathcal{N}}$. Let $F_{\widetilde{\mathcal{N}}}$ be a $\C^{(1)}$-smooth mapping on $\widetilde{\mathcal{N}}$ such that $F_{\widetilde{\mathcal{N}}}|_C = F_{\mathcal{N}}|_C$ around $\Bar{x}$. Define $S = F_{\widetilde{\mathcal{N}}}^{-1} (D)$. Then $S$ is amenable in $\widetilde{\mathcal{N}}$ around $\Bar{x}$. Similarly as before, we have $S = \mathrm{cl} (S_0)$ around $\Bar{x}$, where $S_0 := F_{\widetilde{\mathcal{N}}}^{-1} (\mathrm{int}(D))$. 
    
    Note that $C_0 \subseteq \mathcal{N} \cap \mathcal{M} \subseteq \widetilde{\mathcal{N}}$. Hence $C_0 \cap V \subseteq S_0 \cap V$ for some open neighborhood of $\Bar{x}$ in $\R^n$. Take $V$ to be small enough. We know from Lemma \ref{lem:CQ-manifold} that $S_0 \cap V$ is connected. So $C_0 \cap V = S_0 \cap V$, from where
    \begin{equation*}
        C \cap V = \mathrm{cl} (C_0) \cap V = \mathrm{cl} (S_0) \cap V = S \cap V,
    \end{equation*}
    which means $C$ is also an amenable set in $\widetilde{\mathcal{N}}$ around $\Bar{x}$. Since $\widetilde{\mathcal{N}} \subseteq \mathcal{M}$, by Remark \ref{rmk:amenable-structure} we know $C$ is amenable at $\Bar{x}$ in $\mathcal{M}$.
\finpf

\section{Appendix: prox-regularity on Riemannian manifolds} \label{sec:equi}

The goal of this section is to discuss prox-regularity in connection to other geometric notions that measure how close a set is to being convex, which also allows us to derive the proof of Lemma~\ref{lem:prox-regularity}. In particular, we focus on a characterization of prox-regularity by comparing the intrinsic distance with the Riemannian distance as below (see \cite{haantjes}). 

\begin{defn}
    Given a metric space $(X, d)$ and a set $C \subseteq X$,
    we say that $C$ has {\em finite extrinsic curvature} at a point $\Bar{x} \in C$ in $X$ if there exist $\sigma, r > 0$ such that
    \begin{equation}
        d_C (p, q) \leq d (p, q) + \sigma d (p, q)^3 \qquad \text{for all}~ p, q \in C \cap  B_r (\Bar{x}).
    \end{equation}
\end{defn}

\begin{prop}\label{prop:equiv}
    In a $\C^{(3)}$-smooth Riemannian manifold $\M$, a set $C \subseteq \mathcal{M}$ is prox-regular at a point $\Bar{x} \in C$ if and only if it has finite extrinsic curvature there.
\end{prop}

\pf
    The reverse implication can be deduced from \cite[Remark 3.2]{lewis-lopez-nicolae-22} and \cite[Theorem 1.6]{Lytchak05}. The direct implication follows from Proposition \ref{prop:prox-ast}.
\finpf

\begin{rem}
   Lemma \ref{lem:prox-regularity} is an immediate consequence of the above result. To see this, consider the setting where $\mathcal{M}$ is a $\C^{(3)}$-smooth embedded submanifold of $\R^n$. Using again the inequality (\ref{eq:dm}), we can see that $C$ has finite extrinsic curvature at $\Bar{x}$ in $\R^n$ if and only if it has finite extrinsic curvature at $\Bar{x}$ in $\mathcal{M}$. Lemma \ref{lem:prox-regularity} then follows easily.
\end{rem}

Although equivalences among related properties have been extensively studied in the literature in the setting of Euclidean spaces and Riemannian manifolds, we include here, based on similar techniques, a complete and direct proof of the fact that prox-regularity implies finite extrinsic curvature  with the aim of clarifying details for the reader's convenience. Before doing this, we recall some properties of geodesic metric spaces and smooth Riemannian manifolds.

For $\kappa \in \R$, let $D_\kappa$ be the diameter of the complete, simply connected, $2$-dimensional Riemannian manifold $M_\kappa^2$ of constant sectional curvature $\kappa$. More precisely, $D_\kappa = \infty$ if $\kappa \le 0$, while $D_\kappa = \pi/\sqrt{\kappa}$ if $\kappa > 0$. 

Given three points $x, y, z$ in a geodesic space, a {\em geodesic triangle} $\Delta = \Delta(x,y,z)$ is the union of three geodesic segments (its sides) joining each pair of points. A {\em $\kappa$-comparison triangle} for $\Delta$ is a triangle in $M_\kappa^2$ that has side lengths equal to those of $\Delta$.  

A metric space is called a {\it $\CAT(\kappa)$ space} (also known as a space of curvature bounded above by $\kappa$ in the sense of Alexandrov) if every two points at distance less than $D_\kappa$ can be joined by a geodesic and all geodesic triangles having perimeter less than $2D_\kappa$ are not thicker than the $\kappa$-comparison triangles (in the sense that distances between points on the sides of a geodesic triangle are no larger than corresponding distances in a $\kappa$-comparison triangle). Correspondingly, {\it $\CBB(\kappa)$ spaces} (or spaces that have {\it curvature bounded below by $\kappa$} in the sense of Alexandrov) are defined in a similar way using the condition that all geodesic triangles are not thinner than the $\kappa$-comparison triangles.

For a $\C^{(3)}$-smooth Riemannian manifold $\mathcal{M}$ and a point $\bar x \in \M$, there exists a ball centered at $\Bar{x}$ that is both a $\CAT(\kappa)$ and a $\CBB(\kappa')$ space for some suitable $\kappa, \kappa' \in \R$ (see, e.g., \cite[Chapter~II.1, Appendix]{bridson}, \cite{burago-gromov-per}, \cite{Berest}).

Henceforth, we make the following assumption. Note that the balls and projections below are defined with respect to the Riemannian distance.
\begin{assume} \label{assume-radius}
Consider a $\C^{(3)}$-smooth Riemannian manifold $\mathcal{M}$ and a subset $C$ that is prox-regular at a point $\Bar{x} \in C$. Let $R > 0$ be small enough such that the following conditions are simultaneously met:
\begin{itemize}
    \item $\overline{B}_R (\Bar{x})$ is compact;
    \item $P_C (x)$ is a singleton for all $x \in \overline{B}_R (\Bar{x})$;
    \item $B_R (\Bar{x})$ is $\CAT(\kappa)$, $\CBB(\kappa')$ for some $\kappa > 0, \kappa' < 0$;
    \item $R < \frac{\pi}{2\sqrt{\kappa}}$, so $B_R (\Bar{x})$ is convex and there is a unique geodesic segment connecting every pair of points in $B_R (\Bar{x})$;
    \item for all distinct $x, y \in B_R (\Bar{x})$ with $d_\mathcal{M} (x, y) < R$, the geodesic from $x$ to $y$ can be extended beyond $y$ to a geodesic of length $R$.
\end{itemize}
\end{assume}

We first prove a result on the extension of geodesics in $\mathcal{M}$ (see also \cite{federer, Kleinjohann81}).

\begin{prop}
\label{prop:extension}
   For any $r \in (0, \frac{R}{6}]$, taking arbitrary $x \in \overline{B}_r (\Bar{x}) \backslash C$ and denoting $u = P_C (x)$, we have that the geodesic from $u$ to $x$ can be extended beyond $x$ to a geodesic of length $2 r$, and $u$ is the nearest point in $C$ of all points belonging to this extension.
\end{prop}

\pf We prove the desired property following similar ideas used in \cite{ArizaRuiz16, Vlasov67}. The proof is based on the following claims.
\setcounter{claim}{0}
\begin{claim} \label{cl:ext:1}  
   The mapping $P_C |_{B_{3r} (\Bar{x})}$ is continuous.
\end{claim}
\noindent{\bf Proof of Claim \ref{cl:ext:1}\ \ }
		Suppose $P_C$ is not continuous at some $x \in B_{3r} (\Bar{x})$. Then there exists $\epsilon > 0$ such that for each $k \in \N$ there is a point $x_k \in B_{3r} (\Bar{x})$ with $d_\mathcal{M} (x_k, x) < \frac{1}{k}$ and $d_\mathcal{M} \left( P_C (x_k), P_C (x) \right) \geq \epsilon$. For sufficiently large $k$ such that $\frac{1}{k} < 3r - d_\mathcal{M} (x, \Bar{x})$, we have
            \begin{align*}
                d_\mathcal{M} \left(P_C (x_k), \Bar{x}\right) &\leq d_\mathcal{M} (x_k, C) + d_\mathcal{M} (x_k, \Bar{x}) \\
                &\leq d_\mathcal{M} (x, C) + d_\mathcal{M} (x_k, x) + d_\mathcal{M} (x_k, \Bar{x}) \\
                &\leq d_\mathcal{M} (x, \Bar{x}) + \frac{1}{k} + 3r < 6r \leq R.
            \end{align*}
            Suppose $P_C (x_k) \rightarrow u \in C$. Then
            \begin{equation*}
                d_\mathcal{M} (u, x) = \lim_{k \rightarrow \infty} d_\mathcal{M} \left( P_C (x_k), x_k \right) = \lim_{k \rightarrow \infty} d_\mathcal{M} (x_k, C) = d_\mathcal{M} (x, C).
            \end{equation*}
            Hence, $u = P_C (x)$, which contradicts $d_\mathcal{M} (P_C (x_k), P_C (x)) \geq \epsilon$.
\finpf 
\begin{claim} \label{cl:ext:2}  
   Let $x \in B_{3r} (\Bar{x}) \backslash C$. Denote $u = P_C (x)$. Extend the geodesic from $u$ to $x$ beyond $x$ and, for small $t > 0$, take $x_t$ on this extension such that $d_\mathcal{M} (x_t, u) = (1 + t) d_\mathcal{M} (x, u)$. Then
        \begin{equation*}
            \lim_{t \searrow 0} \frac{d_\mathcal{M} (x_t, C) - d_\mathcal{M} (x, C)}{d_\mathcal{M} (x_t, x)} = 1.
        \end{equation*}
\end{claim}
\noindent{\bf Proof of Claim \ref{cl:ext:2}\ \ } 
	Denote $u_t = P_C (x_t)$. By Claim \ref{cl:ext:1}, $u_t \rightarrow u$. Then
\[d_\mathcal{M} (x_t, u_t) \leq d_\mathcal{M} (x_t, u) \leq d_\mathcal{M} (x_t, x) + d_\mathcal{M} (x, u),\]
from where
\[\frac{d_\mathcal{M} (x_t, u_t) - d_\mathcal{M} (x, u)}{d_\mathcal{M} (x_t, x)} \leq 1.\]
In addition,
\[d_\mathcal{M} (u, x) \leq d_\mathcal{M} (u_t, x) \leq \frac{t}{1 + t} d_\mathcal{M} (u_t, u) + \frac{1}{1 + t} d_\mathcal{M} (u_t, x_t),\]
so $d_\mathcal{M} (x_t, u_t) - d_\mathcal{M} (x, u) \geq t \left(d_\mathcal{M} (x, u) - d_\mathcal{M} (u_t, u) \right).$ This yields
\[\frac{d_\mathcal{M} (x_t, u_t) - d_\mathcal{M} (x, u)}{d_\mathcal{M} (x_t, x)} \geq \frac{d_\mathcal{M} (x, u) - d_\mathcal{M} (u_t, u)}{d_\mathcal{M} (x, u)} = 1 - \frac{d_\mathcal{M} (u_t, u)}{d_\mathcal{M} (x, u)},\]
where $d_\mathcal{M} (u_t, u) \rightarrow 0$.
\finpf
\begin{claim} \label{cl:ext:3}  
   Let $x \in \overline{B}_r (\Bar{x}) \backslash C$ and $\epsilon \in (0, 1]$. Denote $\delta = d_\mathcal{M} (x, C)$. Then there exists $y \in \mathcal{M}$ such that $d_\mathcal{M} (x, y) = 2r - (1 - \epsilon) \delta$ and $d_\mathcal{M} (y, C) \geq 2 r$.
\end{claim}
\noindent{\bf Proof of Claim \ref{cl:ext:3}\ \ } 
            Let $\sigma = \frac{2r - (1 - \epsilon) \delta}{2r - \delta}$ and
            \begin{equation*}
                U = \left\{ y \in \overline{B}_{2r - (1 - \epsilon) \delta} (x) : d_\mathcal{M} (y, x) \leq \sigma \left( d_\mathcal{M} (y, C) - d_\mathcal{M} (x ,C) \right) \right\}.
            \end{equation*}
            Note $U \neq \emptyset$ since $x \in U$. Moreover,
            \begin{equation*}
                d_\mathcal{M} (y, C) \leq d_\mathcal{M} (y, x) + \delta \le 3 r \qquad \text{for all}~y \in \overline{B}_{2r - (1 - \epsilon) \delta} (x).
            \end{equation*}
            Define a partial order $\preceq$ on $U$ by $y_1 \preceq y_2$ if and only if $d_\mathcal{M} (y_1, y_2) \leq \sigma \left( d_\mathcal{M} (y_2, C) - d_\mathcal{M} (y_1 ,C) \right)$.
            Let $(y_i)_{i \in I}$ be a chain in $U$. Then $d_\mathcal{M} (y_i, C) \leq d_\mathcal{M} (y_j, C)$ whenever $i \leq j$, therefore $\left( d_\mathcal{M} (y_i, C) \right)$ converges. So $(y_i)$ is a Cauchy net. Suppose $y_i \rightarrow y_0$. Then $y_0 \in U$ and $y_i \preceq y_0$ for all $i \in I$. Hence, $y_0$ is an upper bound of $(y_i)$. By Zorn's lemma, $U$ has a maximal element $y \in U$. 
            
            Suppose $d_\mathcal{M} (y, x) < 2r - (1 - \epsilon) \delta$. Then
            \begin{equation*}
                d_\mathcal{M} (y, \Bar{x}) \leq d_\mathcal{M} (y, x) + d_\mathcal{M} (x, \Bar{x}) < 3r.
            \end{equation*}
            By Claim \ref{cl:ext:2}, there exists $y_t \neq y$ arbitrarily close to $y$ such that
            \begin{equation*}
                \frac{d_\mathcal{M} (y_t, C) - d_\mathcal{M} (y, C)}{d_\mathcal{M} (y_t, y)} > \frac{1}{\sigma},
            \end{equation*}
            which means $y \preceq y_t$. Then we have $x \preceq y \preceq y_t$ and $y_t \in U$, which leads to contradiction. Hence, $d_\mathcal{M} (y, x) = 2r - (1 - \epsilon) \delta$. Moreover, since $y \in U$ we have
            \begin{equation*}
                d_\mathcal{M} (y, C) \geq \frac{1}{\sigma} d_\mathcal{M} (y, x) + d_\mathcal{M} (x, C) = 2 r.
            \end{equation*}
\finpf
    
    We return to the proof of Proposition \ref{prop:extension}. Let $x \in \overline{B}_r (\Bar{x}) \backslash C$, $u = P_C (x)$ and denote $\delta = d_\mathcal{M} (x, C)$. By Claim \ref{cl:ext:3}, for every $k \in \N$ there exists $y_k \in \mathcal{M}$ with $d_\mathcal{M} (x, y_k) = 2 r - \left(1 - \frac{1}{k}\right) \delta$ and $d_\mathcal{M} (y_k, C) \geq 2 r$. Then
    \begin{equation*}
        2 r \leq d_\mathcal{M} (y_k, C) \leq d_\mathcal{M} (y_k, u) \leq d_\mathcal{M} (y_k, x) + d_\mathcal{M} (x, u) = 2 r + \frac{1}{k} \delta.
    \end{equation*}
    So $d_\mathcal{M} (y_k, C) \rightarrow 2 r$, $d_\mathcal{M} (y_k, u) \rightarrow 2 r$, $d_\mathcal{M} (y_k, x) \rightarrow 2 r - \delta$ and
    \begin{equation*}
        d_\mathcal{M} (y_k, \Bar{x}) \leq d_\mathcal{M} (y_k, x) + d_\mathcal{M} (x, \Bar{x}) < 3r < R. 
    \end{equation*}
    Extend the geodesic from $u$ to $x$ beyond $x$ to a geodesic of length $2 r$ and take $y$ on this extension such that $d_\mathcal{M} (u, y) = 2 r$. Fix any limit point $y'$ of $(y_k)$. We have
    \begin{equation*}
        d_\mathcal{M} (y', C) = d_\mathcal{M}(y', u) = d_\mathcal{M}(y, u) \quad \text{and} \quad d_\mathcal{M}(y', x) = d_\mathcal{M}(y, x).
    \end{equation*}
    By the uniqueness of geodesics and the lower curvature bound, we know $y' = y$. Hence, $d_\mathcal{M} (y, C) = d_\mathcal{M} (y, u)$, which means $P_C (y) = u$.
\finpf

\begin{prop}\label{prop:prox-ast}
    For all sufficiently small $r > 0$, there exists a constant $\sigma > 0$ such that
    \begin{equation*}
        d_{C \cap \overline{B}_r (\Bar{x})} (p, q) \leq d_\mathcal{M} (p, q) + \sigma d_\mathcal{M} (p, q)^3 \qquad \text{for all}~ p, q \in C \cap \overline{B}_r (\Bar{x}).
    \end{equation*}
\end{prop}

\pf
	Suppose that $R > 0$ satisfies, in addition to Assumption \ref{assume-radius}, the following conditions:
\begin{itemize}
    \item for all $x \in B_R (\Bar{x})$, $d_{B_R (\Bar{x})} (\cdot, x)^2$ is $(-\frac{5}{6})$-convex;
    \item for all $t \in [0, \sqrt{-\kappa'} R]$, we have 
    \[\frac{3}{4} t \leq \tanh (t) \leq t, \quad t \leq \sinh (t) \leq \frac{5}{4} t \quad \text{and} \quad 1+\frac{1}{2}t^2 \le \cosh(t) \le  1+\frac{1}{2}t^2+\frac{1}{8}t^4.\]
\end{itemize}
Take $r \leq \min \{\frac{R}{12},\frac{1}{2\sqrt{-\kappa'}}\}$. The proof is based on the following claims.

\setcounter{claim}{0}
\begin{claim} \label{cl:prox:1}  
Let $p, q \in C\cap \overline{B}_r (\Bar{x})$, $m$ be the midpoint of $p$ and $q$, and $u = P_C (m)$. Then $u \in C \cap \overline{B}_r (\Bar{x})$ and
    \begin{equation*}
        d_\mathcal{M} (m, u) \leq \frac{5}{48 r} d_\mathcal{M} (p, q)^2.
    \end{equation*}
\end{claim}
\noindent{\bf Proof of Claim \ref{cl:prox:1}\ \ } We can suppose $m 
\notin C$. By the $(-\frac{5}{6})$-convexity of $d_{B_R (\Bar{x})} (\cdot, \Bar{x})^2$,
    \begin{equation}\label{eq:m-xbar-2}
        d_\mathcal{M} (m, \Bar{x})^2 \leq \frac{1}{2} d_\mathcal{M} (p, \Bar{x})^2 + \frac{1}{2} d_\mathcal{M} (q, \Bar{x})^2 - \frac{5}{24} d_\mathcal{M} (p, q)^2 \leq r^2 - \frac{5}{24} d_\mathcal{M} (p, q)^2.
    \end{equation}
    Applying Proposition \ref{prop:extension}, the geodesic from $u$ to $m$ can be extended beyond $m$ to a geodesic of length $4r$ and $u$ is the nearest point in $C$ of all points belonging to this extension. Take $x$ on this extension such that $d_\mathcal{M} (u, x) = 4r$. 

The sum of the adjacent angles $\angle_{m} (x, p)$ and $\angle_{m} (x, q)$ is $\pi$. Without loss of generality assume $\angle_{m} (x, p) \leq \frac{\pi}{2}$. Denote 
\[a=\sqrt{-\kappa'}d_\mathcal{M} (m, u),\quad b=\sqrt{-\kappa'}d_\mathcal{M} (m, x),\quad c=\sqrt{-\kappa'}d_\mathcal{M} (m, p),\quad \text{and} \quad s=\sqrt{-\kappa'}d_\mathcal{M} (u, x).\]
Then, by the cosine law in $M_{\kappa'}^2$,
    \begin{align*}
        \cosh(s) \le \cosh(\sqrt{-\kappa'}d_\mathcal{M} (x, p)) &\le \cosh(b) \cosh(c) = \cosh(s - a) \cosh(c) \\
        & = \left( \cosh(s) \cosh(a) - \sinh(s) \sinh(a) \right) \cosh(c) \\
        & \leq \left( \cosh(s) \cosh(c) - \sinh(s) \sinh(a) \right) \cosh(c).
    \end{align*}
    We get that $\sinh(s) \sinh(a) \cosh(c) \leq \cosh(s)\sinh(c)^2$, from where 
    \[\frac{3}{4}s\cdot a\le \tanh(s) \sinh(a) \leq \tanh(c) \sinh(c)\le \frac{5}{4}c^2.\]
 	Thus,
    \begin{equation*}
        d_\mathcal{M} (m, u) \leq \frac{5 d_\mathcal{M} (m, p)^2}{3 d_\mathcal{M} (u, x)} = \frac{5}{48 r} d_\mathcal{M} (p, q)^2.
    \end{equation*}
    Hence, by (\ref{eq:m-xbar-2}),
    \begin{equation*}
        d_\mathcal{M} (u, \Bar{x}) \leq d_\mathcal{M} (m, \Bar{x}) + d_\mathcal{M} (m, u) \leq \sqrt{r^2 - \frac{5}{24} d_\mathcal{M} (p, q)^2} + \frac{5}{48 r} d_\mathcal{M} (p, q)^2 < r,
    \end{equation*}
    which means $u \in C \cap \overline{B}_r (\Bar{x})$.
\finpf

\begin{claim}\label{cl:prox:2}      
    Let $p,q,u \in \overline{B}_r (\Bar{x})$. Suppose that the midpoint $m$ of $p$ and $q$ satisfies
    \begin{equation*}
        d_\mathcal{M} (m, u) \leq K d_\mathcal{M} (p, q)^2
    \end{equation*}
    for some constant $K > 0$. Then
    \begin{equation*}
        d_\mathcal{M} (p, u) + d_\mathcal{M} (q, u) \leq d_\mathcal{M} (p, q) + K' d_\mathcal{M} (p, q)^3,
    \end{equation*}
    where $K' = 5 K^2 - \kappa'/32$.
\end{claim}
\noindent{\bf Proof of Claim \ref{cl:prox:2}\ \ } 
Denote 
\[a=\sqrt{-\kappa'}d_\mathcal{M} (p, u),\quad b=\sqrt{-\kappa'}d_\mathcal{M} (q, u),\quad c=\sqrt{-\kappa'}d_\mathcal{M} (m, p),\quad \text{and} \quad s=\sqrt{-\kappa'}d_\mathcal{M} (u, m).\]
Applying the cosine law in the geodesic triangles $\Delta(p,m,u)$ and $\Delta(q,m,u)$, we obtain    
 \[\cosh(a) + \cosh(b) \le 2 \cosh(c) \cosh(s).\]
Using now the bounds on $\cosh(t)$ for $t \in [0, \sqrt{-\kappa'} R]$, 
    \begin{align*}
        2 + \frac{1}{2}\left( a^2+b^2\right) &\leq 2 \left( 1 + \frac{1}{2} c^2 + \frac{1}{8} c^4 \right)  \left( 1 + \frac{1}{2} s^2 + \frac{1}{8} s^4 \right) \\
        &\leq 2 \left( 1 + \frac{1}{2} c^2 + \frac{1}{8} c^4 \right) \left( 1 + s^2 \right) \quad \text{since } s \le  2\sqrt{-\kappa'}r \le 2\\
        & \leq 2 + c^2 + \frac{1}{4} c^4 + 2s^2 + 2c^2s^2 \quad \text{since } c \le  \sqrt{-\kappa'}r \le 2.
    \end{align*}
    Therefore, $a^2+b^2 \le 2c^2 + \frac{1}{2} c^4 + 4s^2 + 4c^2s^2$. Since $c = \frac{\sqrt{-\kappa'}}{2}d_\mathcal{M} (p, q)$ and $s \le \sqrt{-\kappa'}Kd_\mathcal{M} (p, q)^2$, we get that
    \begin{align*}
        d_{\mathcal{M}} (p, u)^2 + d_{\mathcal{M}} (q,u)^2 & \leq \frac{1}{2} d_{\mathcal{M}} (p, q)^2 + \left(4 K^2 - \frac{\kappa'}{32} \right) d_{\mathcal{M}} (p, q)^4 - \kappa' K^2 d_{\mathcal{M}} (p, q)^6 \\
        &\leq \frac{1}{2} d_{\mathcal{M}} (p, q)^2 + \left(5 K^2 - \frac{\kappa'}{32} \right) d_{\mathcal{M}} (p, q)^4, 
    \end{align*}
    where the last inequality follows since $(- \kappa') d_{\mathcal{M}} (p, q)^2 \le  4(-\kappa')r^2 \le 1$. Consequently,
    \begin{align*}
        d_{\mathcal{M}} (p, u) + d_{\mathcal{M}} (q,u) &\leq \sqrt{2 d_{\mathcal{M}} (p, u)^2 + 2 d_{\mathcal{M}} (q,u)^2} \\
        & \leq \sqrt{d_{\mathcal{M}} (p, q)^2 + \left(10 K^2 - \frac{\kappa'}{16} \right) d_{\mathcal{M}} (p, q)^4} \\
        & \leq d_{\mathcal{M}} (p, q) + \left(5 K^2 - \frac{\kappa'}{32} \right) d_{\mathcal{M}} (p,q)^3.
    \end{align*}
    This completes the proof of Claim \ref{cl:prox:2}.
\finpf

Now we continue with the proof of Proposition \ref{prop:prox-ast}.
    Fix any $p, q \in C \cap \overline{B}_r (\Bar{x})$. We will construct a continuous curve in $C \cap \overline{B}_r (\Bar{x})$ whose length satisfies the desired inequality. 
    
    For any points $x, y \in B_R (\Bar{x})$, denote the geodesic in $\mathcal{M}$ from $x$ to $y$ as $\gamma_{x, y}$ and its length as $l_{x,y}$. We recursively define a sequence of continuous curves $\gamma_n : [0, 1] \rightarrow \mathcal{M}$, where $n \in \N \cup \{0\}$, as follows:
    \begin{itemize}
        \item Define $\gamma_0$ as $\gamma_0 (t) = \gamma_{p, q} (t \cdot l_{p,q} )$ for all $t \in [0, 1]$.
        
        \item For $n \geq 0$, denote $p^n_k := \gamma_n (k \cdot 2^{-n})$ for $0 \leq k \leq 2^n$ and consider the points $u_k^n = P_C(m^n_k)$, where $m^n_k$ is the midpoint of $p^n_{k-1}$  and $p^n_k$ and $1 \leq k \leq 2^n$.
        Then define $\gamma_{n+1}$ as
        \begin{equation*}
            \gamma_{n+1} (t) = 
            \begin{cases}
                \gamma_{p^n_{k-1}, u^n_k} \big( (t -\frac{k-1}{2^n}) \cdot 2^{n+1}l_{p^n_{k-1}, u^n_k} \big), & t \in [\frac{k-1}{2^n}, \frac{2k-1}{2^{n+1}}] \\[6pt]
                \gamma_{u^n_k, p^n_k} \big( (t -\frac{2k-1}{2^{n+1}}) \cdot 2^{n+1}l_{u^n_k, p^n_k} \big), & t \in [\frac{2k-1}{2^{n+1}}, \frac{k}{2^n}]
            \end{cases}
        \end{equation*}
        for $1 \leq k \leq 2^n$.
    \end{itemize}
    Using Claim \ref{cl:prox:1} and induction, we can see that for all $n \geq 0$ and $1 \leq k \leq 2^n$,
    \begin{equation}\label{eq:gamma-c}
        \gamma_n \big( \frac{k}{2^{n}} \big) \in C \cap \overline{B}_r (\Bar{x})
    \end{equation}
    and
    \begin{equation*}
        d_\mathcal{M} (m^n_k, u^n_k) \leq K l_{p^n_{k-1}, p^n_k}^2,
    \end{equation*}
    where $K = \frac{5}{48r}$. Thus,
    \begin{equation}\label{eq:m-u}
        \max\{l_{p^n_{k-1}, u^n_k},~ l_{u^n_k, p^n_k}\} \leq \frac{1}{2} l_{p^n_{k-1}, p^n_k} + K l_{p^n_{k-1}, p^n_k}^2.
    \end{equation}
    Moreover, applying Claim \ref{cl:prox:2}, we obtain
    \begin{equation}\label{eq:peri}
        l_{p^n_{k-1}, u^n_k} + l_{u^n_k, p^n_k} \leq l_{p^n_{k-1}, p^n_k} + K' l_{p^n_{k-1}, p^n_k}^3,
    \end{equation}
    where $K' = 5 K^2 - \frac{\kappa'}{32}$. We will use (\ref{eq:m-u}) and (\ref{eq:peri}) to show that $\{\gamma_n\}$ converges to the desired curve in $C \cap \overline{B}_r (\Bar{x})$. 
    
    First, it follows easily from (\ref{eq:m-u}) that
    \begin{equation*}
        \max\Big\{l_{p^{n+1}_{k-1}, p^{n+1}_k}: 1 \leq k \leq 2^{n+1}\Big\} \leq \frac{3}{4} \max\Big\{l_{p^n_{k-1}, p^n_k}: 1 \leq k \leq 2^n\Big\}\qquad \text{for all}~ n \geq 0,
    \end{equation*}
    and hence
    \begin{equation}\label{eq:l}
        l_{p^n_{k-1}, p^n_k} \leq l_{p,q} \big(\frac{3}{4}\big)^n \leq 2r \big(\frac{3}{4}\big)^n \qquad \text{for all}~ 1 \leq k \leq 2^n,~ n \geq 0.
    \end{equation}
    Note that
    \begin{equation}\label{eq:n+1-n}
        \gamma_{n+1} \big( \frac{k}{2^{n}} \big) = \gamma_n \big( \frac{k}{2^{n}} \big) \qquad \text{for all}~0 \leq k \leq 2^n,~ n \geq 0.
    \end{equation}
    Fix any $n \geq 0$ and $1 \leq k \leq 2^n$. For all $t \in [\frac{k-1}{2^n}, \frac{k}{2^n}]$, using (\ref{eq:l}) and (\ref{eq:n+1-n}), we obtain
    \begin{align*}
        d_\mathcal{M} \big( \gamma_{n+1} (t), \gamma_n (t) \big)
        &\leq d_\mathcal{M} \Big( \gamma_{n+1} (t), \gamma_{n+1} \big( \frac{k}{2^{n}} \big) \Big) + d_\mathcal{M} \Big( \gamma_n (t), \gamma_n \big( \frac{k}{2^{n}} \big) \Big) \\
        &\leq l_{p^n_{k-1}, u^n_k} + l_{u^n_k, p^n_k} + l_{p^n_{k-1}, p^n_k} \leq 5r \big(\frac{3}{4}\big)^n.
    \end{align*}
    Therefore, given $n \geq 0$,
    \begin{equation*}
        \max_{0 \leq t \leq 1} d_\mathcal{M} \big( \gamma_{n+1} (t), \gamma_n (t) \big) \leq 5r \big(\frac{3}{4}\big)^n,
    \end{equation*}
    which means that $\{\gamma_n\}$ is a uniformly Cauchy sequence and therefore converges uniformly to a continuous curve $\gamma: [0, 1] \rightarrow \mathcal{M}$. Moreover, it is easy to see from (\ref{eq:gamma-c}) that $\gamma(t) \in C \cap \overline{B}_r (\Bar{x})$ for all $t \in [0, 1]$.
    
    Second, by (\ref{eq:peri}) and (\ref{eq:l}), we have that for all $n \geq 0$,
    \begin{align*}
        L(\gamma_{n+1})
        &\leq \sum_{k = 1}^{2^n} \left(l_{p^n_{k-1}, p^n_k} + K' l_{p^n_{k-1}, p^n_k}^3\right) \\
        &\leq \sum_{k = 1}^{2^n} l_{p^n_{k-1}, p^n_k} \left(1 + K' l_{p,q}^2 \big(\frac{3}{4}\big)^{2n} \right) \\
        &= \left(1 + K' l_{p,q}^2  \big(\frac{9}{16}\big)^n \right) L(\gamma_n) \leq e^{K' l_{p,q}^2  (\frac{9}{16})^n} L(\gamma_n),
    \end{align*}
    where $L(\gamma_n)$ denotes the length of $\gamma_n$. Hence,
    \begin{equation*}
        L(\gamma_n) \leq e^{K' l_{p,q}^2 \sum_{k=0}^{n-1} (\frac{9}{16})^k} l_{p,q} \leq e^{\frac{16}{7} K' l_{p,q}^2} l_{p,q} \leq \big(1 + 4 K' l_{p,q}^2 \big) l_{p,q},
    \end{equation*}
    where the last inequality follows from the fact that
    \begin{equation*}
        \frac{16}{7} K' l_{p,q}^2 \leq \frac{16}{7} \big( \frac{125}{48^2 r^2} - \frac{\kappa'}{32} \big) \cdot 4 r^2 < \frac{1}{2} - \frac{2}{7} \kappa' r^2 \leq \frac{4}{7}
    \end{equation*}
    and that $e^t \leq 1 + \frac{7}{4} t$ for all $t \in [0, \frac{4}{7}]$. Letting $n \rightarrow \infty$, we get
    \begin{equation*}
        d_{C \cap \overline{B}_r (\Bar{x})} (p, q) \leq L(\gamma) \leq l_{p,q} + 4 K' l_{p,q}^3.
    \end{equation*}
    We conclude that the constant $\sigma = 4 K'$ satisfies the required inequality.
\finpf

We finish by mentioning another two related properties defined for sets on manifolds, namely ``2-convexity" and ``positive reach". A closed set $S \subseteq \mathcal{M}$ is 2-convex (see \cite{Lytchak05}) if there exist $\sigma, \rho > 0$ such that
\begin{equation*}
    d_S (p, q) \leq d_\mathcal{M} (p, q) + \sigma d_\mathcal{M} (p, q)^3 \qquad \text{for all}~ p, q \in S ~\text{with}~ d_\mathcal{M} (p, q) < \rho,
\end{equation*} 
and $S$ has positive reach (see \cite{federer}) if there exists $\delta > 0$ such that $P_S (y)$ is a singleton for all $y \in \mathcal{M}$ with $d_\mathcal{M} (y, S) < \delta$. For sets of positive reach, a result related to Proposition \ref{prop:prox-ast} is \cite[Proposition~2.4]{Lytchak24}.

\begin{prop}
    The following statements are equivalent:
    \begin{enumerate}
        \item $C$ is prox-regular at $\Bar{x}$;
        \item $C$ has finite extrinsic curvature at $\Bar{x}$;
        \item there exists $r > 0$ such that $C \cap \overline{B}_r (\Bar{x})$ is 2-convex;
         \item there exists $r > 0$ such that $C \cap \overline{B}_r (\Bar{x})$ has positive reach.
    \end{enumerate}
\end{prop}

\pf
    It follows directly from the definitions that $3 \Rightarrow 2$ and $4 \Rightarrow 1$. Moreover, by \cite[Remark~3.2]{lewis-lopez-nicolae-22}, we know $2 \Rightarrow 3$;  by \cite[Theorem 1.6]{Lytchak05}, we know $3 \Rightarrow 4$. Proposition \ref{prop:equiv} gives the equivalence between $1$ and $2$. We conclude that the four properties are equivalent.
\finpf

\end{document}